\documentclass[USenglish]{article}
\usepackage{fullpage}
\usepackage{amsmath,amsfonts,amsthm,amssymb}
\usepackage{color,xcolor}
\usepackage{ifpdf}
\usepackage{psfrag}
\usepackage{graphicx,graphics}
\usepackage{hhtensor}
\usepackage{comment}
\usepackage[small]{caption}
\usepackage{subfigure}
\usepackage{tikz}
\usepackage{mathtools}
\usepackage{bbm}
\usepackage{algorithm,algorithmic}
\usepackage{pgfplots}
\usepackage{pgfplotstable}
\pgfplotsset{compat = 1.3}

\usepackage{xparse}
\usepackage{bigints}

\usepackage{url}
\usepackage{hyperref}

\usepackage[capitalise,noabbrev]{cleveref}
\usepackage{todonotes}
\newtheorem{theorem}{Theorem}[section]

\newtheorem{proposition}[theorem]{Proposition}
\newtheorem{remark}[theorem]{Remark}

\usepackage{cancel}
\usepackage{algorithm,algorithmic}
\usepackage{enumitem}

\newcommand{\diff}[1]{{\mathrm{d}{#1}}}
\newcommand{\bu}{\mathbf{u}}
\newcommand{\bFF}{\mathcal{F}}
\newcommand{\bGG}{\mathcal{G}}
\newcommand{\bU}{\mathbf{U}}
\newcommand{\bH}{\mathbf{H}}

\newcommand{\LL}{\mathcal{L}}
\newcommand{\xip}{x_{i+1/2}}
\newcommand{\xin}{x_{i-1/2}}
\newcommand{\iip}{{i+1/2}}
\newcommand{\iin}{{i-1/2}}
\newcommand{\bS}{\mathcal{S}}
\usepackage{bbm}

\def\R{\mathbb{R}}

\definecolor{darkspringgreen}{rgb}{0., 0.55, 0.3}
\definecolor{dartmouthgreen}{rgb}{0.05, 0.5, 0.06}
\definecolor{etonblue}{rgb}{0.59, 0.78, 0.64}
\definecolor{airforceblue}{rgb}{0., 0.4, 0.66}
\definecolor{arylideyellow}{rgb}{0.91, 0.84, 0.42}
\definecolor{emerald}{rgb}{0.31, 0.78, 0.47}
\definecolor{uclagold}{rgb}{1.0, 0.7, 0.0}
\definecolor{cadmiumorange}{rgb}{0.93, 0.53, 0.18}

\newsavebox{\DelimiterBox}
\newlength{\DelimiterHeight}
\newlength{\DelimiterDepth}
\newsavebox{\ArgumentBox}
\newlength{\ArgumentHeight}
\newlength{\ArgumentDepth}
\newlength{\ResizedDelimiterHeight}
\newlength{\ResizedDelimiterDepth}

\ifluatex

\else

\fi

\newcommand{\bby}{\mathbf{y}}

\hypersetup{
pdftitle={}
pdfauthor={},
pdfpagemode=UseOutlines,
linkbordercolor=0 0 0,
linkcolor=red,
citecolor=blue,
colorlinks=true
}

\begin{document}
\title{Arbitrary High Order WENO Finite Volume Scheme with Flux Globalization for Moving Equilibria Preservation}

\author{Mirco Ciallella$^{(1)}$\thanks{Corresponding author \href{mailto:mirco.ciallella@ensam.eu}{mirco.ciallella@ensam.eu}}, Davide Torlo$^{(2)}$\thanks{\href{mailto:davide.torlo@sissa.it}{davide.torlo@sissa.it}} and Mario Ricchiuto$^{(3)}\thanks{\href{mailto:mario.ricchiuto@inria.fr}{mario.ricchiuto@inria.fr}}$
\\ 
{\small (1): \'Ecole Nationale Sup\'erieure d'Arts et M\'etiers, Institut de M\'ecanique et d'Ing\'enierie I2M, 33400 Talence, France}\\
{\small (2): SISSA mathLab, Mathematics Area, SISSA, via Bonomea 265, Trieste 34136, Italy}\\
{\small (3): INRIA, Univ. Bordeaux, CNRS, Bordeaux INP, IMB, UMR 5251, France}}
\date{\today}
\maketitle

\begin{abstract}

In the context of preserving stationary states, e.g.\ lake at rest and moving equilibria, a new formulation of the shallow water system,
called Flux Globalization has been introduced by Cheng \textit{et al.} (2019). 
This approach consists in including the integral of the source term in the \textit{global flux} 
and reconstructing the new \textit{global flux} rather than the conservative variables. 
The resulting scheme is able to preserve a large family of smooth and discontinuous steady state moving equilibria. 
In this work, we focus on an arbitrary high order WENO Finite Volume (FV) generalization of the global flux approach.
The most delicate aspect of the algorithm is the appropriate definition of the \emph{source flux} (integral of the source term) and 
the quadrature strategy used to match it with the WENO reconstruction of the hyperbolic flux.
When this construction is correctly done, one can show that the resulting WENO FV scheme admits exact discrete steady states characterized by constant global fluxes.
We also show that, by an appropriate quadrature strategy for the source, we can embed exactly some particular steady states, e.g.\ the lake at rest for the shallow water equations.
It can be shown that an exact approximation of global fluxes leads to a scheme with better convergence properties and improved solutions.
The novel method has been tested and validated on classical cases: subcritical, supercritical and transcritical flows.

\end{abstract}
\vspace*{2mm}


\section{Introduction}
The   Saint-Venant or shallow water equations model  the dynamics of hydrostatic free surface waves under the action of gravity   \cite{de1871theorie}. 
This system of non-linear partial differential equations is  valid under the hypothesis of very large wavelengths, or very shallow depth, and describes the evolution of the 
water depth and volume flux.  They are used in a variety of  engineering applications going from river and estuarine hydrodynamics, to urban flood and tsunami risk assessment.
The basic set of equations constitutes a system of hyperbolic conservation laws, endowed with all the usual properties \cite{Serre99}:  real eigenstructure of the flux Jacobians, 
Rankine-Hugoniot relations defining weak solutions, and  a convex entropy extension to detect the admissible ones. 
The influence of spatial changes of the height of the bottom topography, of the viscous effects in boundary layers, as well as other
physical effects, are accounted for by means of appropriately defined source terms. The addition of these terms results in a system of balance laws which may admit quite 
a large number of  equilibrium solutions  dictated by the interaction between different components of the flux variation in space, and the forcing terms. 
Analytical solutions are not always available for such states, but they can often be characterized in some implicit form (see e.g.  \cite{delestre2013swashes}).\\

%
%

The numerical approximation of the shallow water equations is a very active research domain. It is quite impossible to mention all the relevant literature,   numerous original methods having been   devised in many   different settings: 
finite volumes  \cite{audusse2004fast,gallardo2007well,noelle2007high,diaz2013high,bollermann2013well,kurganov2002central,kurganov2018finite,cheng2016moving,xing2011advantage,berthon2016fully,michel2017well,michel2021two,ciallella2021arbitrary}, finite elements \cite{https://doi.org/10.1002/fld.1650211009,Hauke:1998,https://doi.org/10.1002/fld.250,cite-keyi2,https://doi.org/10.1002/fld.1153,cite-key1,Scovazzi3,BEHZADI2020112662}, 
residual distribution  \cite{ricchiuto2009stabilized,ricchiuto2007application,ricchiuto2015explicit,ricchiuto2011c,arpaia2018r,ArR:20}, discontinuous Galerkin \cite{xing2014exactly,mantri2021well,arpaia2022efficient,meister2014unconditionally} and so on. 
As in many other fields, the final goal is to design robust numerical methods, providing  physically relevant solutions in realistic applications with the lowest possible computational effort.
In this respect, very high order methods are known for their computational efficiency, especially  when dealing with sufficiently smooth solutions.
Among these, an excellent compromise between high order of accuracy and robustness for    non-smooth solutions or data is provided by
WENO finite volume methods   \cite{shu1998essentially, shu1988efficient} which are the object of the present paper.\\


For a given accuracy order, and on a given mesh, a rewarding strategy to design numerical methods with reduced error is 
that of the so called property-preserving discretizations. These methods try to mimic additional consistency conditions than just those
expressed by the system of equations themselves. A typical example of such derived conditions is the energy/entropy (in)equality.
For the shallow water equations, the presence of the source terms makes this issue particularly interesting. As mentioned, the system of balance  laws
admits a large number of equilibrium solutions, depending on the interaction flux-variation/source. This notion is related to the notion of a well-balanced discretization,
usually defined as one capable of reproducing one or more of these equilibria at the discrete level.
The most classical version of this notion is related to the  discrete conservation of  the so called lake at rest state 
with flat free surface and zero velocity \cite{bv94,VAZQUEZCENDON1999497,doi:10.1137/0733001}.
Since then, this issue is essential to the design of any numerical method for   the shallow water equations. 

The numerical evidence shows that, besides  being exact for the lake at rest state, well balanced methods also provide
error reduction for other states. This has pushed the community toward the development of a more general notion of well balancedness,
hoping to further improve the cost efficiency of the resulting schemes. These generalizations try to account for other equilibria than just the lake at rest state,
as well as for more physical effects, and source terms.

%
There have been several approaches proposed to this end. 
Among these we mention  the  fully well balanced generalized Riemann solver 
proposed by Berthon and Chalons \cite{berthon2016fully,berthon2022very}, the optimization based  corrected reconstruction method by Castro and Par\'es \cite{castro-pares-jsc20,GOMEZBUENO2021125820,gomez2021collocation},
the steady residual subtraction technique proposed in 
\cite{berberich2021high}. 
In the previously mentioned methods only local information is used 
to obtain discrete equilibria. 
Another approach consists of a global flux method first introduced in \cite{gascon2001construction}, 
then applied to design non-linear shock capturing methods 
for balance laws in \cite{donat2011hybrid,caselles2009flux}, and more
recently as an alternative way to design fully well balanced methods 
\cite{chertock2018well,cheng2019new,chertock2022well}.
The main idea of these techniques is to build a quasi-conservative hyperbolic system 
starting from a balance law by integrating the source term so that it can be included 
in a global definition of the flux.   
Some methods require the explicit knowledge of the exact equilibrium, some others are agnostic of the specific form of the equilibrium, but might need to solve nonlinear equations at each cell and time, or might not preserve the steady states in all configurations, e.g.\ shocks.
In this work we focus on the flux globalization framework, which belong to the latter class.  This approach relies on 
the definition of a single flux (a global flux) whose divergence expresses the effects of both the divergence of the conservative flux, and that of the forcing term.
By standard techniques it is possible to devise a scheme which is exactly consistent with a constant global flux. This gives a generalized  consistency 
condition which naturally embeds the effects of the source integral. In its simplest formulation, this approach is a natural way to embed in the finite volume setting
the notion of   residual distribution, as discussed in \cite{Abgrall2022}.\\
%


In this work we propose a possible strategy  to combine  flux globalization with very high order of accuracy via a   WENO  approach.
The heart of the method is based on three  main elements:  a high accurate quadrature allowing to define 
cell averages for the global flux;  a  global flux WENO reconstructions based on these averages; a standard upwind numerical flux.
We can easily show that this strategy guarantees the exact consistency with constant cell-averages of the global flux. 
For the Saint-Venant equations with the bathymetric and  friction terms, we consider 
further improvements     to guarantee the exact preservation of the lake at rest solution.  Although the method is not exactly well balanced with respect to all 
analytically computable equilibria, we show that for a wide range of solutions  the overall gains in accuracy obtained with the proposed method are
of several orders of magnitude on a given mesh.\\

The paper is organized as follows. In \cref{se:SW_Equation} we introduce the shallow water equations and the steady state equilibria that we want to preserve. In \cref{sec_Space} we introduce the spatial discretization including the finite volume discretization, the global flux definition, the lake-at-rest well-balanced modification and the high order WENO reconstruction. In \cref{se_time_discretization} we describe the temporal discretization, based on a    Deferred Correction  method. 
In \cref{se:numerics} we perform several tests showing the superior capability of the global flux high order method with respect to other high order methods.
Finally, in \cref{se:summary} we draw some conclusions and we describe some perspective future works.


\section{Shallow Water Equations}\label{se:SW_Equation}
\subsection{Model}

\begin{figure}
	\centering
	\begin{tikzpicture}
		
		\draw [-stealth] (-0.2,-0.3) -- (4.2,-0.3);
		\node [black,scale=1] at (4.2,-0.5) {$x$};
		
		\draw [stealth-stealth] (1,-0.3) -- (1,0.45);
		\node [black,scale=0.9] at (1.5,0) {$b(x)$};
		\draw [stealth-stealth] (3,0.1) -- (3,2.1);
		\node [black,scale=0.9] at (3.2,1.2) {$h$};
		\draw [stealth-stealth] (3.9,-0.3) -- (3.9,2.15);
		\node [black,scale=0.9] at (4.1,1.0) {$\eta$};
		
		\draw [-latex] (1,1) -- (2,1);
		\node [black,scale=0.9] at (1.5,1.2) {$u$};
		
		\draw [black,ultra thick] plot [smooth] coordinates {(0,0) (0.5,0.2) (1,0.5) (1.5,0.3) (2,0.4) (2.5,0.1) (3,0.05) (3.5,0.2) (4,0.15)};
		\draw [blue, ultra thick]  plot [smooth] coordinates {(0,2) (0.5,2.2) (1,2.7) (1.5,2.1) (2,1.7) (2.5,2.1) (3,2.15) (3.5,2.4) (4,2.15)};
		
	\end{tikzpicture}
	\caption{Shallow Water Equations: definition of the variables.}\label{fig:SWE}
\end{figure}
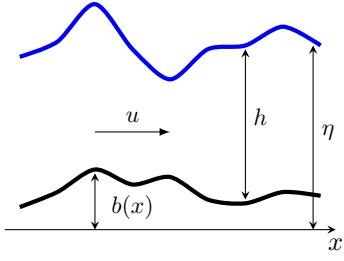

Although the method proposed can be used for more general systems of balance laws, in this paper we will focus on the shallow water equations (SWEs) with
bathymetry and the friction  source terms. For time independent topographies, the system reads 
\begin{equation}\label{eq:CL}
	\frac{\partial \bu}{\partial t} + \frac{\partial\bFF(\bu)}{\partial x} = \bS(\bu,x), \;\;\;\;\;\text{on}\;\;\;\;\;\boldsymbol{\Omega}_T = \boldsymbol{\Omega}\times[0,T]\subset\R\times\R^+,  
\end{equation}
with conserved variables, flux and source term given by
\begin{equation}\label{eq:SWE}
	\bu=\begin{bmatrix} h \\ q  \end{bmatrix}\;,\;\;
	\bFF(\bu)=\begin{bmatrix} q \\ \frac{q^2}{h}+g\frac{h^2}{2} \end{bmatrix}\;,\;\;
	\bS(\bu,x)= \begin{bmatrix} 0 \\ S(\bu,x)  \end{bmatrix} = -gh\begin{bmatrix} 0 \\ \frac{\partial b(x)}{\partial x} \end{bmatrix}
	-gq\begin{bmatrix} 0 \\ \frac{n^2}{h^{7/3}}|q| \end{bmatrix},
\end{equation}
where (cf. Figure~\ref{fig:SWE}) $h$ represents the relative water height, $q$ is its discharge (equal to $hu$, where $u$ is the vertically averaged velocity), 
$g$ is the gravitational acceleration, $b(x)$ is the local bathymetry and $n$ is the Manning friction coefficient. 
It is also convenient to introduce the free surface water level $\eta:=h+b$. \\

The key idea used in this paper is to recast~\eqref{eq:CL}  in the following equivalent pseudo-conservative form:
\begin{equation}\label{eq:globalCL}
	\frac{\partial \bu}{\partial t} + \frac{\partial \bGG(\bu,x)}{\partial x} = 0 \quad\text{ such that }\quad \bGG(\bu,x)=\begin{bmatrix} q \\ K\end{bmatrix}=\begin{bmatrix} q \\ \frac{q^2}{h}+g\frac{h^2}{2}+\mathcal R \end{bmatrix},
\end{equation}
having set 
\begin{equation}\label{eq:globalR}
	\mathcal R(x,t) := -\int^x S(\bu,\xi)\;\diff{\xi}= g \int^x \left[ h(\xi,t)\frac{\partial b(\xi)}{\partial \xi} + \frac{n^2}{h^{7/3}(\xi,t)}|q(\xi,t)|q(\xi,t)\right] \;\diff{\xi}.
\end{equation}
In this case both components of the global flux,  $q$ and $K$, are steady  equilibrium variables. Moreover,
the form of system~\eqref{eq:globalCL}  allows the use of some  classical techniques for conservation laws, being the source term be treated through the definition and the local discrete approximation  of $\mathcal R$.  

\subsection{Steady state equilibria}

%
%


As recalled in the introduction,  system~\eqref{eq:SWE} admits a variety of steady equilibria dependent on the interaction between the variations of  $\bFF$, or of some of its components, and the source.
We are interested in  well balanced  schemes that are capable of preserving in some sense some of these equilibria at the discrete level. 

Classically, the most studied equilibrium is the lake at rest  given by 
\begin{equation}\label{lake_at_rest}
	u=0; \qquad \eta(x,t) =h(x,t) +b(x)  \equiv\eta_0 \in \R, \quad \forall\,x\in \boldsymbol{\Omega},\, t \in [0,T].
\end{equation}
this can also be considered as a special case of a  constant energy state  which is not at rest, defined as 
\begin{equation}\label{nrg}
	h(x,t)u(x,t)\equiv q_0; \qquad \mathcal{E} =g\eta(x,t)   + \dfrac{u^2}{2}  \equiv \mathcal{E}_0 \in \R, \quad \forall\,x\in \boldsymbol{\Omega},\, t \in [0,T].
\end{equation}
There are several works treating specifically these two solutions, and other works deal with other analytical
states, see e.g.  \cite{ricchiuto2015explicit,michel2017well}.

A more general standpoint is to consider  that steady state solutions are generally characterized by the invariants
\begin{equation}\label{eq:newmovequi}
	h(x,t)u(x,t)\equiv q_0 \quad  \text{ and } \quad K(x,t)=\frac{q^2}{h} + g\frac{h^2}{2}+  \int_{x_0}^x g \left[ h(\xi,t)\frac{\partial b(\xi)}{\partial \xi} + \frac{n^2}{h^{7/3}(\xi,t)}|q(\xi,t)|q(\xi,t)\right] \diff \xi \equiv K_0. 
\end{equation}
Obviously, the lake at rest~\eqref{lake_at_rest}  and constant energy~\eqref{nrg} are   a special cases of~\eqref{eq:newmovequi}. For a smooth solution in the frictionless case, for example, we have 
\begin{align}
	q&\equiv q_0,\\
	\begin{split}
		0&=\partial_x \left( \frac{q^2}{h} + g\frac{h^2}{2} \right) + gh \partial_x b   =-\frac{q_0^2}{h^2}\partial_x h + gh\partial_x h+ gh \partial_x b= h \partial_x \left( \frac{q_0^2}{2h^2}+ g (h+b) \right) .
	\end{split}
\end{align}
in which the last reduces to the second in \eqref{nrg}.
Then, the last equation sums up to verify that $\Upsilon:=\frac{q^2}{2h^2}+ g \eta \equiv\Upsilon_0$ is constant in space and time.
However, this condition might not hold when bathymetry or water height are discontinuous in some points, e.g.\ transcritical flows with shock, while \eqref{eq:newmovequi} is always true.

%

\section{Space discretization: Global Flux Finite Volume method}\label{sec_Space}

The global flux idea is discussed in some detail in \cite{cheng2019new} in the context of a piecewise linear approximation. 
In this section we discuss a possible strategy to construct arbitrarily high order extensions based on WENO reconstructions.
To begin with, the computational domain $\boldsymbol{\Omega}$ is discretized
into $N_x$ equispaced control volumes $\Omega_i = [\xin,\xip]$ of size $\Delta x$ centered at $x_i=i\Delta x$ with $i=i_\ell,\ldots,i_r$.\\
For the control volume $\Omega_i$ we can define the cell average at time $t$: 
\begin{equation}
	\bar\bU_i(t):=\frac{1}{\Delta x}\int_{\xin}^{\xip} \bu(x,t)\;\diff{x}.
\end{equation}
The semi-discrete finite volume scheme for the system~\eqref{eq:globalCL} reads
\begin{equation}\label{eq:FV}
	\frac{\diff{\bar\bU_i}}{\diff{t}} + \frac{1}{\Delta x}(\widehat \bH_{\iip}-\widehat \bH_{\iin}) = 0 ,
\end{equation}
where $\widehat \bH_{\iip}$ is a numerical flux consistent with the global flux $\mathcal G$. The global flux differs from the original flux and this makes tricky the development of
an upwind scheme based on the solution of Riemann problems by an approximate solver.  In this work we are thus
focusing on a relatively simple  approach in which upwinding is defined by the homogeneous system. In other words we take
%
\begin{align}\label{eq:upwind flux}
	\widehat \bH_{\iip} &= (L^{-1}\Lambda^+L)_{\iip} \bGG^L_{\iip}  + (L^{-1}\Lambda^-L)_{\iip} \bGG^R_{\iip}.
\end{align}
Here, $\bGG^{L,R}_{\iip}$ are the discontinuous reconstructed point values of the global flux $\mathcal{G}(\bU)$ respectively at the left and right side of the cell interface $x_{\iip}$. $L$ is the matrix of the left eigenvectors computed from the flux Jacobian of the hyperbolic problem~\eqref{eq:CL} in the Roe averaged state, i.e.,
\begin{equation}
	J(\bU^*) = \begin{pmatrix}
		0 & 1\\ -u_*^2 +gh_* & 2u_*
	\end{pmatrix}, \qquad \text{with }\begin{cases}
		h_*=\frac{h^L+h^R}{2},\\
		u_*=\frac{\sqrt{h^L}u^L+\sqrt{h^R}u^R}{\sqrt{h^L}+\sqrt{h^R}}.
	\end{cases}
\end{equation} 
$\Lambda^\pm$ correspond to the upwinding weights
\begin{equation}
	\Lambda_i^+ = \begin{cases} 1,& \text{if } \lambda_i > 0, \\ 0,& \text{if } \lambda_i < 0, \end{cases} \qquad 
	\Lambda_i^- = \begin{cases}1,& \text{if } \lambda_i < 0, \\ 0,& \text{if } \lambda_i > 0. \end{cases}
\end{equation}
It is important to notice that the Jacobian of the flux cannot be directly computed as not all the quantities are available at the interface (we reconstruct only the global flux, not the conserved quantities). Hence, one has to recover the value of $h$ from the global flux itself. To do so, in a positive manner, we use the technique proposed in \cite{cheng2019new} and recalled for completeness  in \cref{sec:positiveHreconstruction}.

Other Riemann solvers, e.g. Rusanov, can be used as long as the numerical dissipation term depending on the difference of the states at the interface vanishes when the steady states are reached. For example, this can be achieved using cutoff functions~\cite{chertock2018well}.

In order to obtain the reconstructed values at the left and right side of interfaces, we will use a high order WENO reconstruction technique
that will be applied directly to cell averages of the global flux. As we will see, this will require
introducing WENO polynomials for 
all the quantities of interests, i.e., $h$, $q$, $b$, $\mathcal{R}$, $\mathcal{G}$ and $K$. 
Note in particular that the simple averaging technique for $\mathcal{R}$ used  in~\cite{cheng2019new}
is not enough to reach more than second order of accuracy.
Another delicate point is  the definition of the local values of the source term used to obtain $\mathcal{R}$. We deal with this aspect
in \cref{sec:larSourceDiscretization}.
\subsection{Global flux evaluation and cell quadrature}
As a starting point for the WENO procedure, we  need cell averages   of the 
global flux, denoted here by
\begin{equation}\label{eq:cellave global flux}
\bar\bGG_i(\bu,x)=\bar\bFF_i(\bu) + \begin{bmatrix} 0 \\ \bar{\mathcal R}_i\end{bmatrix} .
\end{equation}
The first step in our discretization is to build a WENO reconstruction of the conserved variables, hereafter denoted 
by $\tilde{\bu}(x)$, which is used as a basis for the quadrature formula
\begin{equation}\label{eq:flux_average}
\bar{\mathcal{F}}_i(\bu) = \sum_{q} w_q \mathcal{F}(\tilde{\bu}(x_{i,q})),
\end{equation}
where $(x_{i,q},w_q)$ are the high order quadrature points and weights used in the cell $[x_{i-1/2},x_{i+1/2}]$. In particular, we will use a 4-point Gauss--Legendre quadrature.
A similar formula is used for the  integral source $\mathcal{R}$:
\begin{align}\label{eq: Rave}
\bar{\mathcal R}_i \approx \sum_q w_q \mathcal R_{i,q}.
\end{align}
To obtain the values at the quadrature points, we suppose a piecewise polynomial reconstruction of the integral source $\mathcal R$. 
Exploiting the relation   $\partial_x\mathcal{R}\equiv  S$ and we can evaluate the required terms as
\begin{equation}\label{eq: Rquad}
\mathcal R_{i,q} = \mathcal R^R_{\iin} - \int_{\xin^R}^{x_{i,q}} \tilde S(x)\;\diff{x} = \mathcal R^R_{\iin} -  
\sum_\theta \underbrace{\int_{\xin^R}^{x_{i,q}} \ell_\theta(x)\;\diff{x}}_{r_\theta^q} \;S(x_{i,\theta}), \quad i\geq i_l,
\end{equation}
where $\ell_\theta$ are the Lagrangian polynomials associated to  the quadrature points, $r_\theta^q$ are computed exactly,  and having implicitly assumed that
%
\begin{equation}\label{eq:sourcetreatment}
\tilde{S}(x)|_{[x_{i-1/2},x_{i+1/2}]}=\tilde S_i(x) = \sum_\theta \ell_\theta (x) S(x_{i,\theta}).
\end{equation}   
We remark that $q$ and $\theta$ are just free indexes and that both $x_{i,q}$ and $x_{i,\theta}$ are the quadrature nodes on the interval $[x_{i-1/2},x_{i+1/2}]$. Moreover, to maintain the well balancing with respect to lake at rest, a further splitting of the source term has to be done and it will be described in details in Section \ref{sec:larSourceDiscretization}, but, for the moment, we maintain this notation to better convey the main idea of the scheme.

The definition of the array of quadrature values of   $\mathcal{R}$  requires an appropriate  initial value $ \mathcal R^R_{\iin} $.
In practice,  
we impose that $\mathcal{R}$ at the left extrema of the computational domain is equal to 0, i.e.,
\begin{equation}
\mathcal R^R_{i_\ell-1/2} = 0.
\end{equation}
Then, when iterating over the elements to evaluate \eqref{eq: Rquad}  $\forall\, i$,
we set 
\begin{equation}\label{eq:jumpR}
\mathcal R^R_{\iip} = \mathcal R^L_{\iip} + [\![ \mathcal R_{\iip} ]\!],
\end{equation}
where we note that $\mathcal R^L_{\iip}$ is  
\begin{align}\label{eq: Rinterface}
\mathcal R^L_{\iip} = 
\mathcal R^R_{\iin} - \int_{\xin^R}^{\xip^L} \tilde S(x)\;\diff{x} = 
\mathcal R_{\iin}^R - \Delta x \bar S_i,\quad i>1, 
\end{align}
with the average source $\bar{S}_i$ obtained as 
\begin{equation}\label{eq:Save}
\bar S_i := \frac{1}{\Delta x}\int_{\xin^R}^{\xip^L} \tilde S(x)\;\diff{x} = \sum_\theta w_\theta S(x_{i,\theta}).
\end{equation}
Notice that we need to introduce the jump of $\mathcal{R}$ to be able to preserve the lake at rest steady state because the bathymetry is discontinuous at the interfaces. Indeed, for continuous reconstruction of the bathymetry this term vanishes.

Finally, we need to provide a precise definition of how the cell average of the source $\bar{S}_i$ are evaluated, as well as of the jumps $[\![ \mathcal R_{\iip} ]\!]$.
Concerning the latter, for collocated approximations including the boundary nodes, as those used in  \cite{mantri2022}, this jump could be set to zero.
For  the non-collocated  approximation used here this  is not necessarily adequate. For the cell averages $\bar{S}_i$ a special treatment has to be applied to maintain equilibrium for lake at rest simulations. These aspects are discussed  in  \cref{sec:larSourceDiscretization}.
%

\subsection{Weighted Essentially Non-Oscillatory (WENO) reconstruction}\label{sec:WENO}

%
The WENO reconstruction is used in this work for several quantities:
\begin{itemize}
\item for the solution and for the data, in order to be able to define  the fluxes and the source values at quadrature points in \eqref{eq:flux_average}, \eqref{eq: Rquad} and \eqref{eq:Save}
\item for the global flux $\mathcal G(x) = \mathcal F(x) +\mathcal R(x)$ in order to obtain the left and right interface values for the upwind scheme in \eqref{eq:upwind flux}.
\end{itemize}
For each of these quantities, high order approximations are constructed  starting from cell averages.
%
%
We briefly recall the  basics of the WENO reconstruction 
to obtain the polynomial     $u_i(x)$  in cell $i$. We consider  polynomial reconstructions of order $p$, with $p$ odd.  To construct them, we select stencils
of  $p$ cells around cell $i$:
\begin{equation}
\lbrace\Omega_{l_x},\quad l_x = i-r+1, \dots, i+r-1 \rbrace,
\end{equation}
where $2r-1=p$. On each of these stencils, one constructs a high order  polynomial $p^{HO}$  fulfilling  the constraints
\begin{equation}
\frac{1}{\Delta x}\int_{x_{i-j-1/2}}^{x_{i-j+1/2}}p^{HO}(x) \,\diff{x} = u_{i-j}, \qquad j =-r+1,\dots, r-1,
\end{equation}
and $r$ low order polynomials $p_m(x)$, $m=0,\dots ,r-1$, that fulfill
\begin{equation}
\frac{1}{\Delta x}\int_{x_{i-r+j+m-1/2}}^{x_{i-r+j+m+1/2}} p_m(x)\diff{x} = u_{i-j+m}, \qquad j=1,\dots,r .
\end{equation}
The WENO reconstruction aims at combining the low order polynomials in order to obtain the high order reconstruction in case all the low order polynomials are non-oscillatory, while it will prefer the least oscillatory polynomial in case some of these polynomials show oscillations.
As an example, the WENO reconstruction of order 5 (WENO5) will use a stencil of length $p=5$ with $r=3$ low order reconstructions. 
Therefore, the involved cells span from $i-2$ to $i+2$. \\
%
To achieve  high order accuracy  
optimal  linear weights $d_m$  can be defined such that $\sum_m d_m(x) p_m(x) = p^{HO}(x)$ (see e.g.~\cite{jiang1996efficient,balsara2000monotonicity}).
The linear approximation obtained  in this way  however   suffers  from  oscillations and Gibbs phenomena. To remove these 
artifacts, non-linear weights are introduced, defined $\forall\; m = 0,\ldots,r-1$ as 
%
\begin{equation}
\omega_m = \frac{\alpha_m}{\sum^{r-1}_{k=0}\alpha_k} \;,\;\;\; 	\alpha_k = \frac{d_k}{(\beta_k+\epsilon)^2} \;. 
\end{equation}
In the last expression $\epsilon$ is a small number controlling the variations detected by the weights, and    needed to avoid division by zero ($10^{-6}$ in this paper),  while the $\beta_k$ are   smoothness indicators  defined by
\begin{equation}
\beta_k = \sum_{l=1}^{r-1} \int_{\xin}^{\xip} \left(\frac{\diff{}^l}{\diff{x}^l} p_k(x)\right)^2 \Delta x^{2l-1}\diff{x}\;,\;\;\;k = 0,\ldots,r-1.
\end{equation}
The WENO approximation is finally defined as:
\begin{equation}
\tilde{u}(x) = \sum_{m=0}^{r-1} \omega_m(x) p_m(x).
\end{equation}
The values of the optimal weights $d_m$ and the formulae for computing $\beta_m$ can be found  
in~\cite{jiang1996efficient,balsara2000monotonicity} up to $r=6$. In the following, we will test the reconstruction with orders $p=3$ and $p=5$.


\begin{proposition}[Global flux property and steady states]\label{prop:global_flux} Given a state $\lbrace \bar{\mathbf U}_i \rbrace_{i=i_\ell}^{i_r}$ such that $\bar{\mathcal{G}}_i =\mathcal G_0 $ for all $i$, then  
$\lbrace \bar{\mathbf U}_i \rbrace_{i=i_\ell}^{i_r}$ is a steady state of  scheme \eqref{eq:FV} with a WENO reconstruction of the global flux $\mathcal{G}$.
\begin{proof}
	Since every low order polynomial interpolates the constant values $\mathcal{G}_0$, all polynomials are identically the constant $\mathcal{G}_0$. In this case, for the upwind flux  \eqref{eq:upwind flux}, we   trivially  get
	$\widehat \bH_{ i\pm 1/2} = \mathcal G_0$  and, hence, $\frac{d \bar{\mathbf U}_i}{dt } = 0$, i.e., the thesis.
\end{proof}
\end{proposition}

\begin{remark}[ODE integration analogy at equilibrium] 
The exact solution of the global flux 
scheme is not (or not necessarily) the analytical solution
of the PDE. However, there exists a discrete relation of the steady
state equilibrium solution of the global flux method, but it is highly nonlinear 
due to the high order finite volume formulation and the WENO nonlinear reconstruction.
Indeed, at the equilibrium, the PDE \eqref{eq:globalCL} becomes an ODE and an analogy to 
classical ODE solvers may be made considering our discrete solution. 
In practice, due to WENO reconstruction, the presented global flux method can be 
recast into a sort of multistep method, where the global flux equations depend on the left 
and right neighboring cell averages. 

The analysis is simpler for other approximation choices, as shown for example in \cite{mantri2022} where the same ideas have been 
developed in the framework of a DG-SEM discretization using Gauss-Lobatto interpolation points. 
In this case, an explicit equivalence with collocation Runge-Kutta ODE solvers can be shown, and more precisely with the well known 
LobattoIIIA method. In this case, super-convergence results can be applied to precisely characterize the discrete steady states. 
We have not been able so far to provide similar explicit analogies for the method proposed here, 
due to the intricacy of the WENO reconstruction, combined with Gaussian quadrature of the source term.

\end{remark}

\subsection{Source term quadrature for the  shallow water equations} \label{sec:larSourceDiscretization}

For the shallow water equations, we provide here explicit formulas for the evaluation of the source term at quadrature points. We also discuss the issue
of the jumps at the interface. We focus mostly   on the frictionless case. The treatment of friction is  discussed at the end of the paragraph.  \\

We start by defining the WENO   reconstructions of $h$, $\eta$ and $b$ in the quadrature points of the cell $\Omega_i$, denoting them by  $\tilde{h}_{i,q}$, $\tilde{\eta}_{i,q}$ and $\tilde{b}_{i,q}$.
We assume unduly that the same weights are used  for $h$,  and for the bathymetry $b$, so that the consistency with
constants water level is trivially satisfied for the reconstructions $\eta_0\equiv \tilde \eta_{i,q}=\tilde h_{i,q}+\tilde b_{i,q}$.  In  practice,  the weights are directly computed using $\eta$ 
to detect this state.  Then, we obviously have that, at each quadrature point,
\begin{equation}
\tilde{h}_{i,q} = \tilde{\eta}_{i,q} -\tilde{b}_{i,q} .
\end{equation}
Let us also define a Lagrange interpolation of the bathymetry inside the $i$-th cell and its evaluation at the interfaces:
\begin{equation}\label{eq:bath_reconstruction}
\tilde{b}_i(x):= \sum_{q} \ell_q(x) \tilde{b}_{i,q},\qquad \text{and} \qquad b^L_{\iip} = \tilde{b}_i(x_{\iip}), \,\qquad b^R_{\iin} = \tilde{b}_i(x_{\iin}) .
\end{equation}

We start by looking at the bathymetry term source which is the most critical as it contains a derivative in space.
Following \cite{xing2005high} we start by writing this term as  
\begin{equation}\label{eq:sourceWB}
S(\bu,x) =- g h(x) \partial_x b(x) =- g \eta(x) \partial_x b(x) + g\partial_x\left(\frac{b^2(x)}{2}\right).
\end{equation}
To evaluate the integral of the source we now sample the WENO polynomials of $\eta$ and $b$ in the quadrature points.
More precisely, we compute the required values of $\mathcal{R}$ as 
\begin{align}\label{eq: Rquad2}
\mathcal R_{i,q} &= \mathcal R^R_{\iin} - \int_{\xin^R}^{x_{i,q}} S(\bu(x),x)\;\diff{x} \\&=
\mathcal R^R_{\iin} + g\int_{\xin^R}^{x_{i,q}} \eta(x)\partial_x b(x) \;\diff{x} - g\left(\frac{(b_{i,q})^2}{2} - \frac{(b_{\iin}^R)^2}{2}\right) 
\end{align}
where in each quadrature point  $x_{i,q}$ now we set
$$
\partial_x b(x_{i,q}) = \sum\limits_s \ell'_s(x_{i,q}) b(x_{i,s})\;.
$$
Now, using a Lagrange interpolation of the whole source term in the quadrature points, we obtain
\begin{align}\label{eq: Rquad3}
\mathcal R_{i,q}  = \mathcal{R}_{i-1/2}^R+
g  \int_{\xin^R}^{x_{i,q}}\sum_\theta \ell_\theta(x)\;\tilde{\eta}_{i,\theta}\sum_s \ell'_s(x_\theta) \tilde{b}_{i,s}\;\diff{x} - g\left(\frac{(\tilde b_{i,q})^2}{2} - \frac{(b_{\iin}^R)^2}{2}\right).
\end{align}
Similarly, we can define the cell averages \eqref{eq:Save} of the source as 
\begin{equation}
	\bar{S}_i = \frac{1}{\Delta x } \left[ g  \int_{\xin^R}^{\xip^L}\sum_\theta \ell_\theta(x)\;\tilde{\eta}_{i,\theta}\sum_s \ell'_s(x_\theta) \tilde{b}_{i,s}\;\diff{x} - g\left(\frac{(b_{\iip}^L)^2}{2} - \frac{(b_{\iin}^R)^2}{2}\right)\right],
\end{equation}
so that, the left interface terms, introduced in \eqref{eq: Rinterface}, can be defined as
\begin{equation}
\begin{split}
	\mathcal{R}_{i+1/2}^L=\mathcal{R}_{i-1/2}^R -\Delta x\bar{S}_i.
\end{split}
\end{equation}
Finally, we define the jump of $\mathcal R$ across the interfaces. The objective is to obtain  exactness for lake at rest state, i.e., $K$ has to be constant in this equilibrium. For more details see Appendix \ref{app:lar} and Remark \ref{rem:definition_jump_R}.
This leads to
\begin{equation}\label{eq:jump}
[\![ \mathcal R_{\iip} ]\!] := g\frac{\eta_\iip^R +\eta_\iip^L }{2}\left(b_\iip^R  - b_\iip^L\right) - g\left(\frac{(b_\iip^R)^2}{2} - \frac{(b_\iip^L)^2}{2}\right).
\end{equation}
This strategy is similar to classical strategies used in well-balanced path conservative methods  \cite{CASTRO2017131}, 
where one uses a linear/segment path to connect the left and right states when evaluating the integral.
These  definitions allow to easily prove the following property.
\begin{proposition}[Lake at rest preservation] \label{prop:lar}The global flux WENO finite volume scheme with quadrature \eqref{eq: Rquad3} of the bathymetric source,
and with definition   \eqref{eq:jump}  of the source integral jump at the interface   is exactly well balanced for the lake at rest state.
\begin{proof} See appendix \ref{app:lar}.
\end{proof}
\end{proposition}

\begin{remark}[Global Flux with jump of $\mathcal{R}$]
	Being the jump of $\mathcal R$ included in the definition of the global flux $\mathcal{G}$, it does not affect its preservation when moving equilibria are considered, as Proposition \ref{prop:global_flux} is still valid.
\end{remark}

To conclude,  
when adding the friction term we  only account for its contributions to the cell integrals. We thus set in the general case
\begin{equation}\label{eq:global_with_friction}
\mathcal{R}_{i,q} = \mathcal R^R_{\iin} + g\sum_{\theta} \int_{\xin^R}^{x_{i,q}} \ell_\theta(x)\;\diff{x} \left( \tilde{\eta}_{i,\theta} \sum_s \ell'_s(x_{i,\theta})\tilde{b}_{i,s}+ g\frac{\tilde{q}_{i,\theta}|\tilde{q}_{i,\theta}| n^2}{\tilde{h}_{i,\theta}^{7/3}} \right)  - g\left[\frac{(b_{i,\theta})^2}{2} - \frac{(b_{\iin}^R)^2}{2}\right] 
\end{equation}
and
\begin{equation}\label{eq:global_with_friction_inter}
\!\!\!\mathcal{R}^L_{\iip} \!\!= \mathcal R^R_{\iin} + g\sum_{\theta} \int_{\xin^R}^{x_{\iip}}\!\! \ell_\theta(x)\diff{x} \left( \tilde{\eta}_{i,\theta} \sum_s \ell'_s(x_{i,\theta})\tilde{b}_{i,s}+ g\frac{\tilde{q}_{i,\theta}|\tilde{q}_{i,\theta}| n^2}{\tilde{h}_{i,\theta}^{7/3}} \right) \! - g\left[\frac{(b^L_{\iip})^2}{2} - \frac{(b_{\iin}^R)^2}{2}\right].
\end{equation}

In \cref{algo:globalFluxReconstruction} we summarize the steps of the reconstruction of the source integral.

\begin{algorithm}
\begin{algorithmic}
	\STATE $\mathcal{R}_{i_l-1/2}:=0$
	\FOR{$i=i_l,\dots,i_r$}
	\STATE Reconstruct the variables $h$, $\eta$ and $b$ in each quadrature point $\theta$ of the cell, obtaining $\tilde{h}_{i,\theta}$, $\tilde \eta_{i,\theta}$ and $\tilde b_{i,\theta}$ using the same WENO weights (computed for $\eta$) 
	\STATE Reconstruct $q$ in the quadrature points obtaining $\tilde q_{i,\theta}$
	\STATE Define $\mathcal{R}_{i,q}$ as in \eqref{eq:global_with_friction}
	\STATE Define $\mathcal{R}^L_{\iip}$ as in \eqref{eq:global_with_friction_inter}
	\STATE Define $[\![\mathcal{R}_{\iip}]\!]$ as in \eqref{eq:jump}
	\STATE Define $\mathcal{R}^R_{\iip}:=\mathcal{R}^L_{\iip}+[\![\mathcal{R}_{\iip}]\!]$
	\ENDFOR
\end{algorithmic}
\caption{Source integral reconstruction}\label{algo:globalFluxReconstruction}
\end{algorithm}

\section{Time discretization}\label{se_time_discretization}

Time integration can be performed with any time discretization method that fits the method of line approach, e.g. RK schemes. In this paper, we employ a  Deferred Correction (DeC) method because it is a family of one step methods with arbitrarily high order of accuracy.
%
The original DeC formulation was introduced in~\cite{daniel1968iterated}, then developed and studied in its different forms in~\cite{dutt2000dec, minion2003dec,christlieb2010integral, liu2008strong}. A 
slightly different form was presented in~\cite{abgrall2017high} for applications to finite element methods.  In the reference   the DeC is presented as an iterative procedure that involves two operators. The iteration process mimic the Picard--Lindel\"of proof at the discrete level with a fixed-point iterative method. Each iteration aims at gaining  one order of accuracy,
so that the order of accuracy sought can be reached with a finite number of corrections.\\
%
%

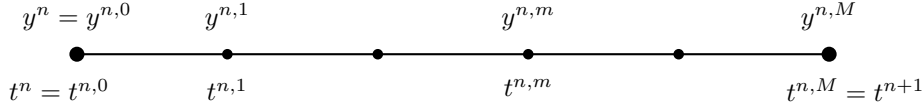
\begin{figure}[h]
	\centering
	\begin{tikzpicture}
		\draw [thick]   (0,0) -- (10,0) node [right=2mm]{};
		\fill[black]    (0,0) circle (1mm) node[below=2mm] {$t^n=t^{n,0} \,\, \quad$} node[above=2mm] {$y^n=y^{n,0}$}
		(2,0) circle (0.7mm) node[below=2mm] {$t^{n,1}$} node[above=2mm] {$y^{n,1}$}
		(4,0) circle (0.7mm) node[below=2mm] {}
		(6,0) circle (0.7mm) node[below=2mm] {$t^{n,m}$} node[above=2mm] {$y^{n,m}$}
		(8,0) circle (0.7mm) node[below=2mm] {}
		(10,0) circle (1mm) node[below=2mm] {$\qquad t^{n,M}=t^{n+1}$} node[above=2mm] {$y^{n,M}$}; 
	\end{tikzpicture} \caption{Time interval divided into sub-time steps}\label{Fig:Time_interval}
\end{figure}

To describe the basic form of the method consider the ODE
\begin{equation}\label{eq:initial_prob}
	\begin{aligned}
		y'(t) = f(y(t)),  \qquad y(t_0)= y_0,
	\end{aligned}
\end{equation}
where $y:\R\to \R^S$ and $f: \R^S\to\R^S$, and in our setting $f$ is  the semidiscretized operator  from \eqref{eq:FV}.  
As usual, the temporal domain is discretized in time steps $[t^n, t^{n+1}]$. We also introduce  $M$ sub-time steps $\lbrace [t^{n,m-1},t^{n,m}]\rbrace_{m=1}^M$ where the boundary points coincide with the extrema of the time step, i.e., $t^n=t^{n,0}$ and $t^{n,M}=t^{n+1}$. The DeC method exploits predictions of the solution  at each  sub-time nodes $t^{n,m}$ denoted as $y^{n,m}$ as explained in Figure~\ref{Fig:Time_interval}.\\

We proceed by first defining a low order approximation of the integral version of the ODE \eqref{eq:initial_prob}, denoted by   $\LL^1$.
The requirement is that the solution of $\LL^1(y)=0$ should be easily obtained. We can for example consider the 
explicit Euler method. We then consider a high order discretization of the ODE, which may be 
costly to solve, as e.g. in high order full tableau implicit RK methods. We denote this operator by $\LL^2$,
and define it as the collocated multi-stage method
\begin{equation}\label{eq:L2}
	\LL^2(y^{n,0}, \dots, y^{n,M}) :=
	\begin{cases}
		y^{n,M}-y^{n,0} - \Delta t \sum_{r=0}^M \theta_r^M f(y^{n,r}),\\
		\vdots\\
		y^{n,1}-y^{n,0} - \Delta t \sum_{r=0}^M \theta_r^1 f(y^{n,r}),
	\end{cases} \approx \begin{cases}
		y^{n,M}-y^{n,0} - \int_{t^{n,0}}^{t^{n,M}}  f(y(s)) \diff s,\\
		\vdots\\
		y^{n,1}-y^{n,0} -  \int_{t^{n,0}}^{t^{n,1}}   f(y(s)) \diff s,
	\end{cases}
\end{equation}
where $\theta_r^m$ is the integral of the $r$-th Lagrangian basis function defined on the sub-time nodes over the interval $[t^{n,0},t^{n,m}]$ \cite{veiga2021dec,abgrall2017high} normalized by the factor $\Delta t$.
This operator is a high order discretization of the integral form of the ODE in each sub-time step. Depending on the chosen sub-time nodes, the order of accuracy varies, for example, with equispaced nodes we can obtain a scheme with order $M+1$, while with Gauss-Lobatto nodes we obtain a $2M$-th accurate scheme, i.e., the Lobatto IIIA schemes \cite{hairer1996solving}. The system of equations $\LL^2=0$ is often a strongly coupled system of nonlinear algebraic equations. This is 
certainly the case when   $f$ is  the nonlinear discretization of a  balance law as the shallow water equations.\\

Concerning the low order operator $\LL^1$, it  consists of an explicit Euler step at each sub-time step:
\begin{equation}\label{eq:L1}
	\LL^1(y^{n,0}, \dots, y^{n,M}) :=
	\begin{cases}
		y^{n,M}-y^{n,0}  - \beta^M \Delta t f(y^{n,0}), \\
		\vdots\\
		y^{n,1}- y^{n,0}- \beta^1 \Delta t f(y^{n,0}),
	\end{cases}
\end{equation}
where $\beta^m= (t^{n,m}-t^{n,0})/\Delta t$.  \\

The DeC iterations provide a new value for the array
of stage values $\bby^{(k)}:=(y^{n,0}, \dots, y^{n,M})^{(k)}$,
given a prediction $\bby^{(k-1)}:=(y^{n,0}, \dots, y^{n,M})^{(k-1)}$.
The method is   first  initiated by setting all entries to the last available solution:  $\bby^{0}=(y(t^n), \dots, y(t^n))$. 
The new iterates are then computed  by evaluating   $\forall\, k=1,\dots,K$  
\begin{equation}
	\label{eq:explicit_dec_explicit}
	\LL^1(\bby^{(k)})=\LL^1(\bby^{(k-1)})-\LL^2(\bby^{(k-1)})
\end{equation}
where $K$ is the final number of iterations. In particular, after $K$ iterations the order of accuracy of the method will be the minimum between $K$ and the accuracy of $\LL^2$. As an example, to obtain order 5, we need $K=5$ and for 
equispaced sub-time nodes $M=4$ (for Gauss--Lobatto nodes $M=3$ suffices). The interested reader, can refer to 
\cite{veiga2021dec,abgrall2021relaxation,torlo2020hyperbolic}  for more details and recent developments.

\section{Numerical Simulations}\label{se:numerics}

The arbitrary high order well-balanced WENO finite volume scheme based on Flux Globalization 
has been tested and validated on several test cases to assess convergence properties and performances.
For comparison, we have used a standard WENO FV scheme with Rusanov numerical flux \cite{shu1998essentially}, without any global flux or other property-preserving features,
with the source term computed using a high order Gauss-Legendre quadrature formula starting from the analytical formulation of the bathymetry and its derivative and the WENO reconstruction of $h$.
We have also investigated the effect of having the global flux property without including the well-balancedness for lake at rest, by computing
the source term with a high order Gauss-Legendre quadrature formula starting from the analytical formulation of the bathymetry and its derivative and the WENO reconstruction of $h$. In this formulation, all the reconstructions are performed with WENO weights obtained by the reconstructed variables themselves, not by the smoothness indicators of the free surface level $\eta$. Still, the global flux is assembled as in \eqref{eq:cellave global flux}.
We do not compare the presented global flux methods with \textit{still}-water well-balanced methods, for example the one
presented in~\cite{shu2006high}. Indeed, as observed in~\cite{xing2011advantage}, such methods do not provide better performances on
moving-water equilibria. 
Both WENO3 and WENO5 schemes are tested. All numerical simulations have been run at CFL $=0.5$.\\

The benchmarks used are quite classical in literature. They allow on one hand to verify first the correct implementation  and convergence 
of the WENO method,  even the non well-balanced. On the other, the results prove 
the dramatic increase in  accuracy obtained with the global flux method for cases at rest, moving, with/out friction,
and for both continuous and discontinuous bathymetries and solutions.
%
%
%

\subsection{Lake at rest}\label{se:lar}
\begin{figure}
	\centering
	\subfigure{
		\begin{tikzpicture}
			\begin{axis}[
				xmin=0,xmax=25,
				grid=major,
				xlabel={$x$},
				ylabel={},
				xlabel shift = 1 pt,
				ylabel shift = 1 pt,
				legend pos= south east,
				legend style={nodes={scale=1, transform shape}},
				tick label style={font=\scriptsize},
				width=.4\textwidth
				]
				\addplot[thick,red]   table [y=eta, x=x]{PERTLARweno5GF_0.dat};
				\addplot[thick,dotted,blue]   table [y expr=\thisrow{b}*0.1, x=x]{PERTLARweno5GF_0.dat};
			\end{axis}
	\end{tikzpicture}}
	\caption{Lake at rest solution: $\eta$ (red) and $b$ (blue).}\label{LAR: solution}
\end{figure}
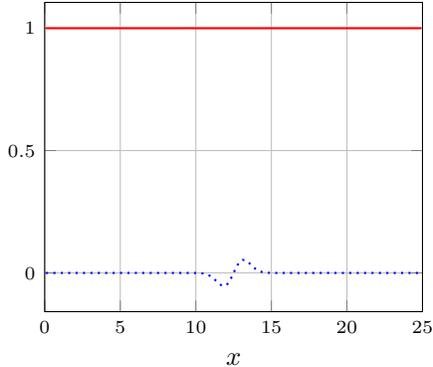
Firstly, we consider the lake-at-rest solution characterized by the initial data
\begin{equation}\label{LAR: IC}
	h(x,0) = 1 - b(x),\qquad q(x,0)\equiv 0,
\end{equation}
over the computational domain $[0,25]$ with subcritical inlet/outlet at the two boundaries. 
The bathymetry employed contains a sinusoidal bump damped at the boundaries (see Figure~\ref{LAR: solution}
for better visualization of solution and bathymetry) and it reads
\begin{equation}\label{LAR: Bath}
	b(x) = 0.05 \sin\left(x-12.5\right)\exp\left(1-(x-12.5)^2\right).
\end{equation}
Let us remark that the classical bathymetry provided in~\cite{delestre2013swashes} is only piecewise polynomial, but globally only $\mathcal{C}^0$. Hence, it is not suited to test the accuracy of very high order methods. On the other side, the bathymetry \eqref{LAR: Bath} is $\mathcal{C}^\infty$ and it has values smaller than machine precision at the boundaries.
\begin{table}
	\caption{Lake at rest: errors and estimated order of accuracy (EOA) with WB and non-WB schemes, using WENO3 and WENO5 reconstructions.}\label{LAR:weno convergence}
	\scriptsize
	\centering
	\begin{tabular}{c|cc|cc|cc|cc} \hline\hline
		&\multicolumn{4}{c|}{Non-WB} &\multicolumn{4}{c}{WB}\\[0.5mm] \hline
		&\multicolumn{2}{c|}{$h$} &\multicolumn{2}{c|}{$q$} &\multicolumn{2}{c|}{$h$} &\multicolumn{2}{c}{$q$}\\[0.5mm]
		\cline{2-9}
		$N_e$ & $L_2$ error        & EOA & $L_2$  error       & EOA & $L_2$ error       &EOA & $L_2$ error      & EOA \\ \hline \hline
		&\multicolumn{4}{c|}{GF-WENO3}&\multicolumn{4}{c}{GF-WENO3}\\ \hline
		25         &  1.0384E-4   &   --        &  4.7943E-5   &  --         &  9.8858E-14  &   --        &  1.2228E-15 &   --    \\
		50         &  1.5496E-5   &    2.67     &  9.2488E-6   &    2.31     &  9.8667E-14  &   --        &  1.4249E-15 &   --    \\
		100        &  1.2117E-6   &    3.62     &  3.6777E-7   &    4.59     &  9.8276E-14  &   --        &  1.6041E-15 &   --    \\
		150        &  2.6776E-7   &    3.69     &  1.5898E-7   &    2.05     &  1.9644E-13  &   --        &  3.3908E-15 &   --    \\
		200        &  9.6323E-8   &    3.53     &  7.6469E-8   &    2.53     &  1.9619E-13  &   --        &  3.6713E-15 &   --    \\
		400        &  8.2671E-9   &    3.53     &  6.0441E-9   &    3.65     &  2.9360E-13  &   --        &  6.1689E-15 &   --    \\
		800        &  6.8811E-10  &    3.58     &  4.7122E-10  &    3.67     &  5.8655E-13  &   --        &  1.3035E-14 &   --    \\ \hline\hline 
		&\multicolumn{4}{c|}{GF-WENO5}&\multicolumn{4}{c}{GF-WENO5}\\  \hline
		25         &  5.1800E-5       &   --        &  6.1657E-5  &  --         &  9.8947E-14  &   --        &  1.3247E-15  &   --    \\
		50         &  4.4066E-6       &   3.45      &  1.5244E-6  &    5.18     &  9.8661E-14  &   --        &  1.4060E-15  &   --    \\
		100        &  6.7998E-7       &   2.66      &  3.5908E-7  &    2.06     &  9.8289E-14  &   --        &  1.5992E-15  &   --    \\
		150        &  1.5437E-7       &   3.63      &  8.8535E-8  &    3.42     &  1.9639E-13  &   --        &  3.4157E-15  &   --    \\
		200        &  4.1973E-8       &   4.50      &  2.3725E-8  &    4.55     &  1.9611E-13  &   --        &  3.7034E-15  &   --    \\
		400        &  1.3952E-9       &   4.89      &  7.5991E-10 &    4.95     &  2.9357E-13  &   --        &  6.2007E-15  &   --    \\
		800        &  4.3120E-11      &   5.01      &  2.2633E-11 &    5.06     &  5.8648E-13  &   --        &  1.3039E-14  &   --    \\
		\hline\hline
	\end{tabular}
\end{table}
The gravitational constant is considered to be $g=1$ and the simulation is run until the final time $T=1$ with
$N_e=\{25,50,100,150,200,400,800\}$ uniform cells using both the non-well-balanced (non-WB) and well-balanced (WB) version of the GF-WENO 
algorithm to assess the convergence and well-balancing properties. The convergence tests performed with WENO3 and WENO5 reconstructions are listed
in Table~\ref{LAR:weno convergence}. It can be noticed that the error decay for the non-WB simulations matches the order 
of the reconstruction for both GF-WENO3 and GF-WENO5, while, for the WB cases, the scheme is able to preserve the exact solution
up to machine precision.

\subsection{Small perturbation of the lake-at-rest solution}\label{se:lar perturbation}

For this test case, we analyze the perturbation of the lake-at-rest solution characterized by  
\begin{equation}\label{LAR perturbation: IC}
	h(x,0) = 1 - b(x) + \begin{cases}  \alpha \psi(x) ,  &\text{if } 9 < x < 10 \\ 0, &\text{otherwise} \end{cases},\qquad q(x,0)\equiv 0
\end{equation}
with $\psi(x)$ a perturbation function defined by
\begin{equation}\label{eq:perturbation}
	\psi(x):= \exp \left( 1- \frac{1}{(1-r(x))^2}\right), \qquad \text{with } r(x) := 4(x-9.5)^2
\end{equation}
and $\alpha =10^{-4}$ over the computational domain $[0,25]$ with subcritical inlet/outlet at the two boundaries.
The bathymetry is a rescaling of~\eqref{LAR: Bath} and it is defined as
\begin{equation}\label{LAR perturbation: Bath}
	b(x) = 0.5 \sin\left(x-12.5\right)\exp\left(1-(x-12.5)^2\right).
\end{equation}
A slightly different bathymetry, with respect to Section~\ref{se:lar}, has been chosen in order to introduce
more noise in the non-WB simulation and appreciate more the method capabilities. 
The simulation was run using a mesh with 150 cells.
The gravitational constant is considered to be $g=9.8$ and the simulation is run until the final time $T=1.5$.\\

\begin{figure}
	\centering
	\subfigure[$t=0$]{
		\begin{tikzpicture}
			\begin{axis}[
				ymin=-1.e-4,ymax=1.e-4,
				xmin=0,xmax=25,
				grid=major,
				xlabel={$x$},
				ylabel={},
				xlabel shift = 1 pt,
				ylabel shift = 1 pt,
				legend pos= south east,
				legend style={nodes={scale=0.6, transform shape}},
				tick label style={font=\scriptsize},
				width=.25\textwidth
				]
				\addplot[red]   table [y=DH, x=x]{PERTLARweno5noGF_N800_0.dat};
				\addplot[thick,dotted,blue]   table [y expr=(\thisrow{b}-1.00000000000)*0.00005, x=x]{PERTLARweno5noGF_N800_0.dat};
			\end{axis}
	\end{tikzpicture}}
	\subfigure[$t=0.5$]{
		\begin{tikzpicture}
			\begin{axis}[
				ymin=-1.e-4,ymax=1.e-4,
				xmin=0,xmax=25,
				grid=major,
				xlabel={$x$},
				ylabel={},
				xlabel shift = 1 pt,
				ylabel shift = 1 pt,
				legend pos= south east,
				legend style={nodes={scale=0.6, transform shape}},
				tick label style={font=\scriptsize},
				width=.25\textwidth
				]
				\addplot[red]   table [y=DH, x=x]{PERTLARweno5noGF_N800_1.dat};
				\addplot[thick,dotted,blue]   table [y expr=(\thisrow{b}-1.00000000000)*0.00005, x=x]{PERTLARweno5noGF_N800_0.dat};
			\end{axis}
	\end{tikzpicture}}
	\subfigure[$t=1$]{
		\begin{tikzpicture}
			\begin{axis}[
				ymin=-1.e-4,ymax=1.e-4,
				xmin=0,xmax=25,
				grid=major,
				xlabel={$x$},
				ylabel={},
				xlabel shift = 1 pt,
				ylabel shift = 1 pt,
				legend pos= south east,
				legend style={nodes={scale=0.6, transform shape}},
				tick label style={font=\scriptsize},
				width=.25\textwidth
				]
				\addplot[red]   table [y=DH, x=x]{PERTLARweno5noGF_N800_2.dat};
				\addplot[thick,dotted,blue]   table [y expr=(\thisrow{b}-1.00000000000)*0.00005, x=x]{PERTLARweno5noGF_N800_0.dat};
			\end{axis}
	\end{tikzpicture}}
	\subfigure[$t=1.5$]{
		\begin{tikzpicture}
			\begin{axis}[
				ymin=-1.e-4,ymax=1.e-4,
				xmin=0,xmax=25,
				grid=major,
				xlabel={$x$},
				ylabel={},
				xlabel shift = 1 pt,
				ylabel shift = 1 pt,
				legend pos= south east,
				legend style={nodes={scale=0.6, transform shape}},
				tick label style={font=\scriptsize},
				width=.25\textwidth
				]
				\addplot[red]   table [y=DH, x=x]{PERTLARweno5noGF_N800_3.dat};
				\addplot[thick,dotted,blue]   table [y expr=(\thisrow{b}-1.00000000000)*0.00005, x=x]{PERTLARweno5noGF_N800_0.dat};
			\end{axis}
	\end{tikzpicture}}
	\caption{Small perturbation of the lake at rest solution computed with the WENO5 scheme: $h-h_{eq}$ (red) and rescaled $b$ (blue) with $N_e=800$.}\label{LAR perturbation_N800:weno5} 
\end{figure}
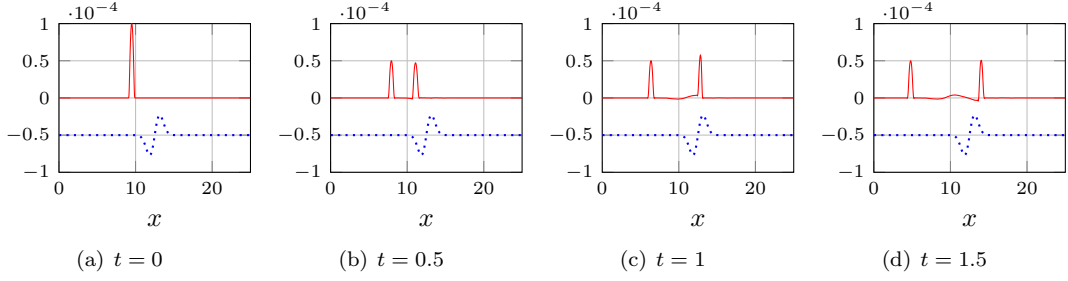

Figure~\ref{LAR perturbation:gf-weno5} shows the evolution over time of the perturbation added over the lake at rest solution
computed with the GF-WENO5 WB scheme.
To better present the results, the plots show the relative variable $h-h_{eq}$ where $h_{eq}$ represents the lake at rest solution
without perturbation provided in Equation~\eqref{LAR: IC}. We then compare the results in Figure~\ref{LAR perturbation:gf-weno5}
with those computed using the classical WENO5 approach, pointed out in Figure~\ref{LAR perturbation:weno5}. It should be noticed
that the classical approach fails at correctly reproducing the perturbation. Indeed, it generates a discretization error higher than the perturbation itself, eventually spoiling the final result. 
Contrary to that, the GF-WENO5 WB scheme correctly reproduce the perturbation which splits into two waves traveling at opposite directions and interacting with the bathymetry. The obtained result is fairly accurate, also considering the coarse mesh used for this study. 
To reach results similar to the GF one with $N_e=150$, we need around $N_e=800$ cells for the classical WENO5 method, see Figure~\ref{LAR perturbation_N800:weno5}.
\begin{figure}
	\centering
	\subfigure[$t=0$]{
		\begin{tikzpicture}
			\begin{axis}[
				ymin=-1.e-4,ymax=1.e-4,
				xmin=0,xmax=25,
				grid=major,
				xlabel={$x$},
				ylabel={},
				xlabel shift = 1 pt,
				ylabel shift = 1 pt,
				legend pos= south east,
				legend style={nodes={scale=0.6, transform shape}},
				tick label style={font=\scriptsize},
				width=.25\textwidth
				]
				\addplot[red]   table [y=DH, x=x]{PERTLARweno5GF_0.dat};
				\addplot[thick,dotted,blue]   table [y expr=(\thisrow{b}-1.00000000000)*0.00005, x=x]{PERTLARweno5GF_0.dat};
			\end{axis}
	\end{tikzpicture}}
	\subfigure[$t=0.5$]{
		\begin{tikzpicture}
			\begin{axis}[
				ymin=-1.e-4,ymax=1.e-4,
				xmin=0,xmax=25,
				grid=major,
				xlabel={$x$},
				ylabel={},
				xlabel shift = 1 pt,
				ylabel shift = 1 pt,
				legend pos= south east,
				legend style={nodes={scale=0.6, transform shape}},
				tick label style={font=\scriptsize},
				width=.25\textwidth
				]
				\addplot[red]   table [y=DH, x=x]{PERTLARweno5GF_1.dat};
				\addplot[thick,dotted,blue]   table [y expr=(\thisrow{b}-1.00000000000)*0.00005, x=x]{PERTLARweno5GF_1.dat};
			\end{axis}
	\end{tikzpicture}}
	\subfigure[$t=1$]{
		\begin{tikzpicture}
			\begin{axis}[
				ymin=-1.e-4,ymax=1.e-4,
				xmin=0,xmax=25,
				grid=major,
				xlabel={$x$},
				ylabel={},
				xlabel shift = 1 pt,
				ylabel shift = 1 pt,
				legend pos= south east,
				legend style={nodes={scale=0.6, transform shape}},
				tick label style={font=\scriptsize},
				width=.25\textwidth
				]
				\addplot[red]   table [y=DH, x=x]{PERTLARweno5GF_2.dat};
				\addplot[thick,dotted,blue]   table [y expr=(\thisrow{b}-1.00000000000)*0.00005, x=x]{PERTLARweno5GF_0.dat};
			\end{axis}
	\end{tikzpicture}}
	\subfigure[$t=1.5$]{
		\begin{tikzpicture}
			\begin{axis}[
				ymin=-1.e-4,ymax=1.e-4,
				xmin=0,xmax=25,
				grid=major,
				xlabel={$x$},
				ylabel={},
				xlabel shift = 1 pt,
				ylabel shift = 1 pt,
				legend pos= south east,
				legend style={nodes={scale=0.6, transform shape}},
				tick label style={font=\scriptsize},
				width=.25\textwidth
				]
				\addplot[red]   table [y=DH, x=x]{PERTLARweno5GF_3.dat};
				\addplot[thick,dotted,blue]   table [y expr=(\thisrow{b}-1.00000000000)*0.00005, x=x]{PERTLARweno5GF_0.dat};
			\end{axis}
	\end{tikzpicture}}
	\caption{Small perturbation of the lake at rest solution computed with the GF-WENO5 WB scheme: $h-h_{eq}$ (red) and rescaled $b$ (blue) with $N_e=150$.}\label{LAR perturbation:gf-weno5}
\end{figure}
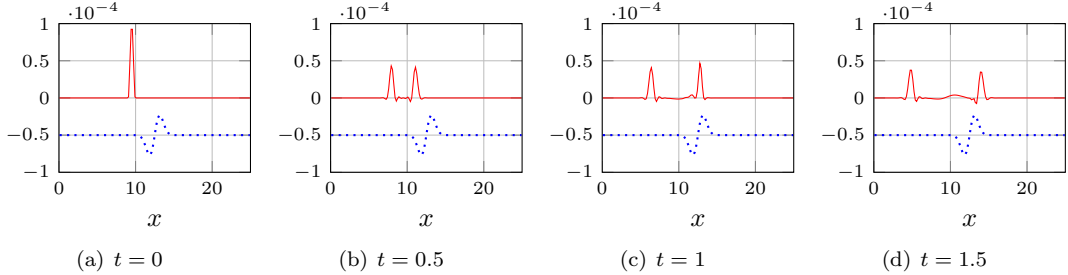
\begin{figure}
	\centering
	\subfigure[$t=0$]{
		\begin{tikzpicture}
			\begin{axis}[
				ymin=-4.e-4,ymax=4.e-4,
				xmin=0,xmax=25,
				grid=major,
				xlabel={$x$},
				ylabel={},
				xlabel shift = 1 pt,
				ylabel shift = 1 pt,
				legend pos= south east,
				legend style={nodes={scale=0.6, transform shape}},
				tick label style={font=\scriptsize},
				width=.25\textwidth
				]
				\addplot[red]   table [y=DH, x=x]{PERTLARweno5noGF_0.dat};
				\addplot[thick,dotted,blue]   table [y expr=(\thisrow{b}-1.00000000000)*0.00020, x=x]{PERTLARweno5GF_0.dat};
			\end{axis}
	\end{tikzpicture}}
	\subfigure[$t=0.5$]{
		\begin{tikzpicture}
			\begin{axis}[
				ymin=-4.e-4,ymax=4.e-4,
				xmin=0,xmax=25,
				grid=major,
				xlabel={$x$},
				ylabel={},
				xlabel shift = 1 pt,
				ylabel shift = 1 pt,
				legend pos= south east,
				legend style={nodes={scale=0.6, transform shape}},
				tick label style={font=\scriptsize},
				width=.25\textwidth
				]
				\addplot[red]   table [y=DH, x=x]{PERTLARweno5noGF_1.dat};
				\addplot[thick,dotted,blue]   table [y expr=(\thisrow{b}-1.00000000000)*0.00020, x=x]{PERTLARweno5GF_0.dat};
			\end{axis}
	\end{tikzpicture}}
	\subfigure[$t=1$]{
		\begin{tikzpicture}
			\begin{axis}[
				ymin=-4.e-4,ymax=4.e-4,
				xmin=0,xmax=25,
				grid=major,
				xlabel={$x$},
				ylabel={},
				xlabel shift = 1 pt,
				ylabel shift = 1 pt,
				legend pos= south east,
				legend style={nodes={scale=0.6, transform shape}},
				tick label style={font=\scriptsize},
				width=.25\textwidth
				]
				\addplot[red]   table [y=DH, x=x]{PERTLARweno5noGF_2.dat};
				\addplot[thick,dotted,blue]   table [y expr=(\thisrow{b}-1.00000000000)*0.00020, x=x]{PERTLARweno5GF_0.dat};
			\end{axis}
	\end{tikzpicture}}
	\subfigure[$t=1.5$]{
		\begin{tikzpicture}
			\begin{axis}[
				ymin=-4.e-4,ymax=4.e-4,
				xmin=0,xmax=25,
				grid=major,
				xlabel={$x$},
				ylabel={},
				xlabel shift = 1 pt,
				ylabel shift = 1 pt,
				legend pos= south east,
				legend style={nodes={scale=0.6, transform shape}},
				tick label style={font=\scriptsize},
				width=.25\textwidth
				]
				\addplot[red]   table [y=DH, x=x]{PERTLARweno5noGF_3.dat};
				\addplot[thick,dotted,blue]   table [y expr=(\thisrow{b}-1.00000000000)*0.00020, x=x]{PERTLARweno5GF_0.dat};
			\end{axis}
	\end{tikzpicture}}
	\caption{Small perturbation of the lake at rest solution computed with the WENO5 scheme: $h-h_{eq}$ (red) and rescaled $b$ (blue) with $N_e=150$.}\label{LAR perturbation:weno5} 
\end{figure}
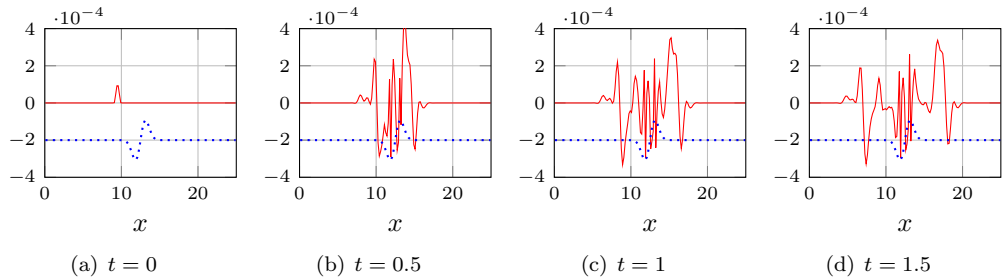

\subsection{Steady states with smooth bathymetry without friction ($n=0$)}\label{se:steadyNoFriction}

Here, we test the method for some moving equilibria steady state problems. We run the tests up to convergence towards steady state in different situations. In the subcritical and supercritical tests the bathymetry is smooth
and equal to~\eqref{LAR: Bath}, as we want to assess the high order accuracy of the schemes. For the transcritical tests we use a modification of the bathymetry used in~\cite{delestre2013swashes}, to study a very similar discontinuous problem.
Depending on the initial and boundary conditions set at the borders of the domain, the flow may be supercritical, 
subcritical or transcritical. 
The meshes used are defined by $N_e$ uniform cells, with different $N_e \in [25,800]$.
We consider the following three sets of final time $T$ and, initial and boundary conditions:
\begin{itemize}
	\item Supercritical flow
	\begin{alignat}{2}\label{eq:super_ic}
		& T = 50,\nonumber\\
		&h(x,0) = 2 - b(x), \qquad &&q(x,0)\equiv 0, \\
		&h(0,t) = 2, 	    \qquad &&q(0,t)=24 \nonumber,
	\end{alignat}
	\item Subcritical flow
	\begin{alignat}{2}\label{eq:sub_ic}
		& T = 200,\nonumber \\
		&h(x,0) = 2 - b(x), \qquad &&q(x,0)\equiv 0, \\
		&q(0,t)=4.42,       \qquad &&h(25,t) = 2\nonumber ,
	\end{alignat}
	\item Transcritical flow
	\begin{equation}
		b(x) = 
		\begin{cases} 
			0.2 \exp\left(1 - \frac{1}{1-\left(\frac{|x-10|}{5}\right)^2} \right), &\text{ if } |x-10|<5, \\ 
			0, &\text{ else}, 
		\end{cases} \label{eq:trans bath}\\
	\end{equation}
	\begin{alignat}{2}\label{eq:trans_ic}
		& T = 200,\nonumber  \\
		&h(x,0) = 0.33 - b(x), \qquad &&q(x,0)\equiv 0, \\
		&q(0,t) = 0.18,        \qquad &&h(25,t)= 0.33.\nonumber
	\end{alignat}
\end{itemize}
The gravitational constant is set to $g=9.812$ for all these tests.
For these three cases, we compare the results obtained using a classical WENO finite volume scheme and the new approach based on flux globalization.
As already mentioned, for supercritical and subcritical cases, the bathymetry~\eqref{LAR: Bath} allows us to perform
convergence tests for very high order methods, when also the flow is smooth.
For the transcritical case with shock, only a qualitative analysis of the test is performed.
Hence, when supercritical and subcritical flows are of interest, we can study the convergence properties of the new scheme by finding the exact solution
given by the non-linear equations taken from~\cite{delestre2013swashes}. Both the WB and non-WB versions of the scheme have
been run to compare the influence of the formulation on the ability of preserving the balanced steady state solution. 
Finally, we also run the same test cases with the
classical WENO3 and WENO5 schemes. \\

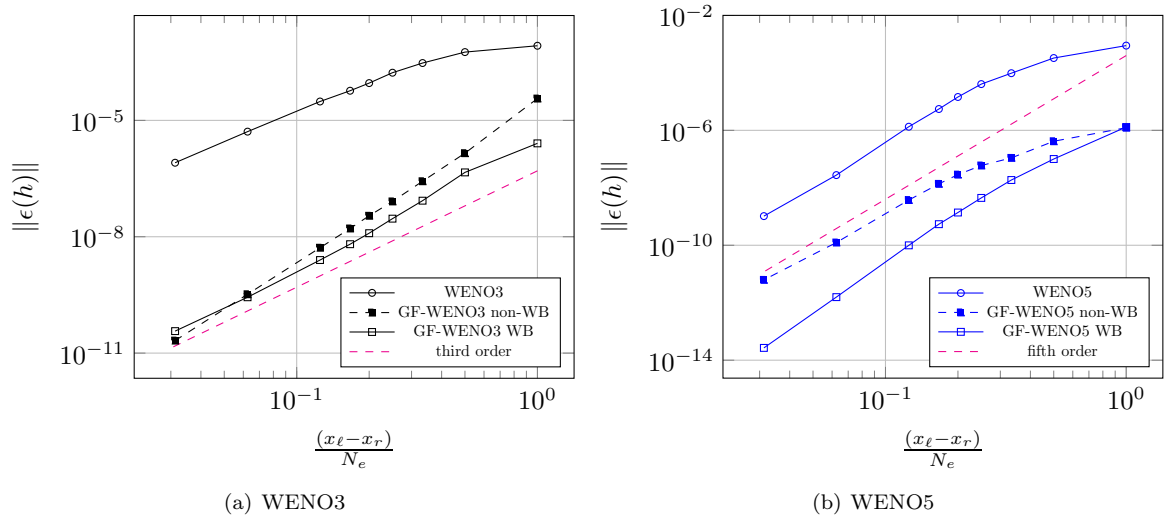
\begin{figure}
	\centering
	\subfigure[WENO3]{
		\begin{tikzpicture}
			\begin{axis}[
				xmode=log, ymode=log,
				grid=major,
				xlabel={$\frac{(x_\ell-x_r)}{N_e}$},
				ylabel={$\|\epsilon(h)\|$},
				xlabel shift = 1 pt,
				ylabel shift = 1 pt,
				legend pos= south east,
				legend style={nodes={scale=0.6, transform shape}},
				width=.45\textwidth
				]
				\addplot[mark=o,mark size=1.3pt,black]             table [y=h   , x=N]{SUPweno3noGF};
				\addlegendentry{WENO3}
				\addplot[mark=square*,dashed,mark size=1.3pt,black]table [y=h   , x=N]{SUPweno3WB};
				\addlegendentry{GF-WENO3 non-WB}
				\addplot[mark=square,mark size=1.3pt,black]        table [y=hWB , x=N]{SUPweno3WB};
				\addlegendentry{GF-WENO3 WB}
				
				\addplot[magenta,domain=1:0.03, dashed]{0.0000005*x*x*x};
				\addlegendentry{third order}
				
			\end{axis}
	\end{tikzpicture}}
	\subfigure[WENO5]{
		\begin{tikzpicture}
			\begin{axis}[
				xmode=log, ymode=log,
				grid=major,
				xlabel={$\frac{(x_\ell-x_r)}{N_e}$},
				ylabel={$\|\epsilon(h)\|$},
				xlabel shift = 1 pt,
				ylabel shift = 1 pt,
				legend pos= south east,
				legend style={nodes={scale=0.6, transform shape}},
				width=.45\textwidth
				]
				\addplot[mark=o,mark size=1.3pt,blue]               table [y=h   , x=N]{SUPweno5noGF};
				\addlegendentry{WENO5}
				\addplot[mark=square*,dashed,mark size=1.3pt,blue]  table [y=h   , x=N]{SUPweno5WB};
				\addlegendentry{GF-WENO5 non-WB}
				\addplot[mark=square,mark size=1.3pt,blue]          table [y=hWB , x=N]{SUPweno5WB};
				\addlegendentry{GF-WENO5 WB}
				
				\addplot[magenta,domain=1:0.03, dashed]{0.0004*x*x*x*x*x};
				\addlegendentry{fifth order}
				
			\end{axis}
	\end{tikzpicture}}
	\caption{Supercritical flow: convergence tests with WENO3 and WENO5.}\label{SUP: convergence}
\end{figure}

\begin{figure}
	\centering
	\subfigure[WENO3]{
		\begin{tikzpicture}
			\begin{axis}[
				xmode=log, ymode=log,
				grid=major,
				xlabel={$\frac{(x_\ell-x_r)}{N_e}$},
				ylabel={$\|\epsilon(h)\|$},
				xlabel shift = 1 pt,
				ylabel shift = 1 pt,
				legend pos= south east,
				legend style={nodes={scale=0.6, transform shape}},
				width=.45\textwidth
				]
				\addplot[mark=o,mark size=1.3pt,black]             table [y=h   , x=N]{SUBweno3noGF};
				\addlegendentry{WENO3}
				\addplot[mark=square*,dashed,mark size=1.3pt,black]table [y=h   , x=N]{SUBweno3WB};
				\addlegendentry{GF-WENO3 non-WB}
				\addplot[mark=square,mark size=1.3pt,black]        table [y=hWB , x=N]{SUBweno3WB};
				\addlegendentry{GF-WENO3 WB}
				
				\addplot[magenta,domain=1:0.06, dashed]{0.00001*x*x*x};
				\addlegendentry{third order}
				
			\end{axis}
	\end{tikzpicture}}
	\subfigure[WENO5]{
		\begin{tikzpicture}
			\begin{axis}[
				xmode=log, ymode=log,
				grid=major,
				xlabel={$\frac{(x_\ell-x_r)}{N_e}$},
				ylabel={$\|\epsilon(h)\|$},
				xlabel shift = 1 pt,
				ylabel shift = 1 pt,
				legend pos= south east,
				legend style={nodes={scale=0.6, transform shape}},
				width=.45\textwidth
				]
				\addplot[mark=o,mark size=1.3pt,blue]               table [y=h   , x=N]{SUBweno5noGF};
				\addlegendentry{WENO5}
				\addplot[mark=square*,dashed,mark size=1.3pt,blue]  table [y=h   , x=N]{SUBweno5WB};
				\addlegendentry{GF-WENO5 non-WB}
				\addplot[mark=square,mark size=1.3pt,blue]          table [y=hWB , x=N]{SUBweno5WB};
				\addlegendentry{GF-WENO5 WB}
				
				\addplot[magenta,domain=1:0.13, dashed]{0.005*x*x*x*x*x};
				\addlegendentry{fifth order}
				
			\end{axis}
	\end{tikzpicture}}
	\caption{Subcritical flow: convergence tests with WENO3 and WENO5.}\label{SUB: convergence}
\end{figure}
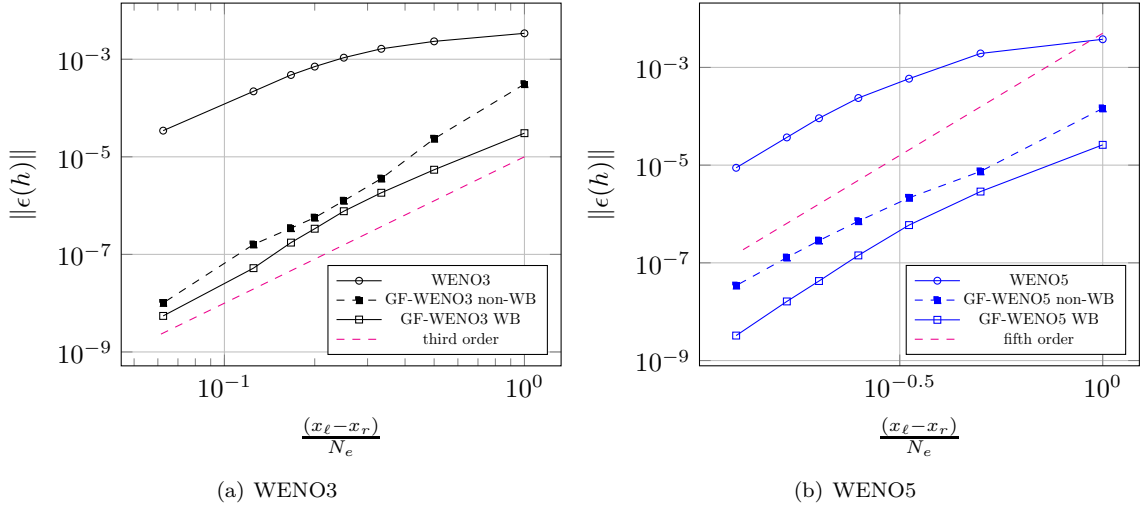

Convergence curves for supercritical and subcritical flows are depicted in Figures~\ref{SUP: convergence} and~\ref{SUB: convergence}, respectively.
All curves, both for WENO3 and WENO5, show the correct third and fifth order accuracy. However, it should be noticed that the GF formulation allows a much better prediction of the solution with errors dropping from 2  to 5 orders of magnitude. 
This is particularly striking for the well balanced global flux  method exactly preserving the lake at rest state.
To have a fairer comparison of the methods, one should compare the computational time with respect to the obtained error. 
Anyway, the methods have not been written in order to optimize the computational costs, so the comparison could be improved. We notice that for a fixed mesh we need more or less double the time of the classical FV scheme to run the GF-WB simulation.  
Let us consider fifth order methods and suppose that the errors of global flux are $10^4$ times smaller than the classical FV ones. Then, the computational time needed for classical FV schemes to obtain the same error of accuracy is 20 times larger than the GF-WB ones.\\

%

\begin{figure}
\centering
\subfigure[$\eta$]{
	\begin{tikzpicture}
		\begin{axis}[
			xmin=0,xmax=25,
			grid=major,
			xlabel={$x$},
			ylabel={},
			xlabel shift = 1 pt,
			ylabel shift = 1 pt,
			legend pos= north east,
			legend style={nodes={scale=0.6, transform shape}},
			tick label style={font=\scriptsize},
			width=.25\textwidth
			]
			\addplot[thick,dotted,blue]  table [y=b   , x=x]{TRANSweno5GF.dat};
			\addplot[red]   table [y=eta   , x=x]{TRANSweno5GF.dat};
			\addplot[black,densely dashed] table [y expr=\thisrow{h}+\thisrow{b}   , x=x]{TRANSweno5noGF.dat};
		\end{axis}
\end{tikzpicture}}
\subfigure[$q$]{
	\begin{tikzpicture}
		\begin{axis}[
			xmin=0,xmax=25,
			grid=major,
			xlabel={$x$},
			ylabel={},
			xlabel shift = 1 pt,
			ylabel shift = 1 pt,
			tick label style={font=\scriptsize},
			ytick distance=0.03,
			y tick label style={/pgf/number format/.cd,fixed relative,precision=6},
			legend pos= south east,
			legend style={nodes={scale=0.6, transform shape}},
			width=.25\textwidth
			]
			\addplot[red]   table [y=q   , x=x]{TRANSweno5GF.dat};
			\addplot[black,densely dashed] table [y expr=\thisrow{h}*\thisrow{u}   , x=x]{TRANSweno5noGF.dat};
		\end{axis}
\end{tikzpicture}}
\subfigure[$\Upsilon$]{
	\begin{tikzpicture}
		\begin{axis}[
			xmin=0,xmax=25,
			grid=major,
			xlabel={$x$},
			ylabel={},
			xlabel shift = 1 pt,
			ylabel shift = 1 pt,
			tick label style={font=\scriptsize},
			ytick distance=0.5,
			y tick label style={/pgf/number format/.cd,fixed relative,precision=6},
			legend pos= south east,
			legend style={nodes={scale=0.6, transform shape}},
			width=.25\textwidth
			]
			\addplot[red]   table [y expr=0.5*\thisrow{q}^2/\thisrow{h}^2+9.812*(\thisrow{h}+\thisrow{b})   , x=x]{TRANSweno5GF.dat};
			\addplot[black,densely dashed] table [y expr=\thisrow{u}^2/2+9.812*(\thisrow{h}+\thisrow{b})   , x=x]{TRANSweno5noGF.dat};
		\end{axis}
\end{tikzpicture}}
\subfigure[$K$]{
	\begin{tikzpicture}
		\begin{axis}[
			xmin=0,xmax=25,
			grid=major,
			xlabel={$x$},
			ylabel={},
			xlabel shift = 1 pt,
			ylabel shift = 1 pt,
			tick label style={font=\scriptsize},
			ytick distance=0.002,
			y tick label style={/pgf/number format/.cd,fixed relative,precision=4},
			legend pos= south east,
			legend style={nodes={scale=0.6, transform shape}},
			width=.25\textwidth
			]
			\addplot[red]   table [y=k, x=x]{TRANSweno5GF.dat};
		\end{axis}
\end{tikzpicture}}
\caption{Transcritical flow: relevant variables computed with GF-WENO5 (red continuous line), WENO5 (black dashed line) schemes and $b$ (blue dotted line) with $N_e=100$.}\label{TRANS: solution}
\end{figure}
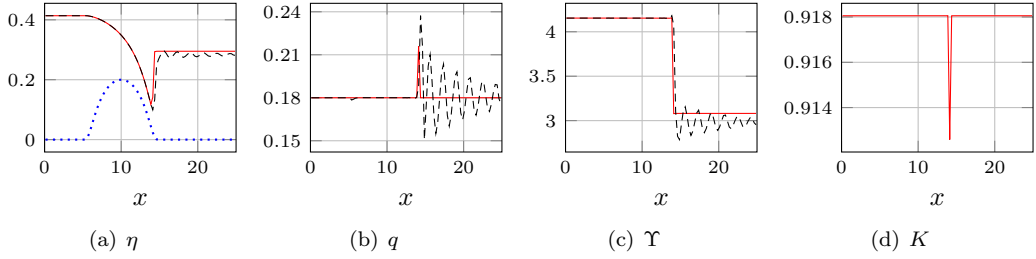

After the convergence analysis, we focus on a more qualitative analysis of the results computed with the WENO5 and GF-WENO5 methods. 
To do so we introduce a set of four representative variables that we are going to use to compare the classical WENO5 with the new GF-WENO5 WB.
The three main variables that we used to assess the new method are $\eta$, $q$ and $K$ ($K$ is only defined for GF formulation). In particular, it should be noticed that, by construction, 
the methods based on flux globalization preserve moving equilibria, i.e., $q_x\equiv 0,\,K_x\equiv 0$. This means that the GF-WENO5 will approximate these quantities up to the order of the residual of the time derivative at the end of the simulation.
In our case, we are able to preserve the constant $q$ and $K$ up to $\sim 10^{-9}$. The last variable we decided to study is the aforementioned $\Upsilon$, which
corresponds to the smooth formulation of the more general global flux introduced herein. Although our approach is developed to preserve other equilibria,
this variable allows to get more insights about the capability of the new algorithm since, for smooth flows, the analytical $\Upsilon$ should be constant at equilibrium.
Solutions for supercritical and subcritical flows are shown in Figures~\ref{SUP: solution} and~\ref{SUB: solution}. For both cases, it
is clear that $q$ and $K$ are well-preserved, and $\Upsilon$ is much better predicted with respect to the one computed through the classical WENO5 method.
The test case that stands out more among the three situations considered here is the transcritical flow with shock in Figure~\ref{TRANS: solution}. 
The WENO5 method
introduces spurious oscillations where the shock occurs that are then propagated in the rest of the computational domain making the overall solution spoiled.
By using the new GF-WENO5 approach, oscillations are not present and the correct solution is recovered before and after the shock.
It can be noticed that, when discontinuous flows are of interest, $\Upsilon$ is not globally constant but features a jump where the shock occurs, see Figure~\ref{TRANS: solution}.
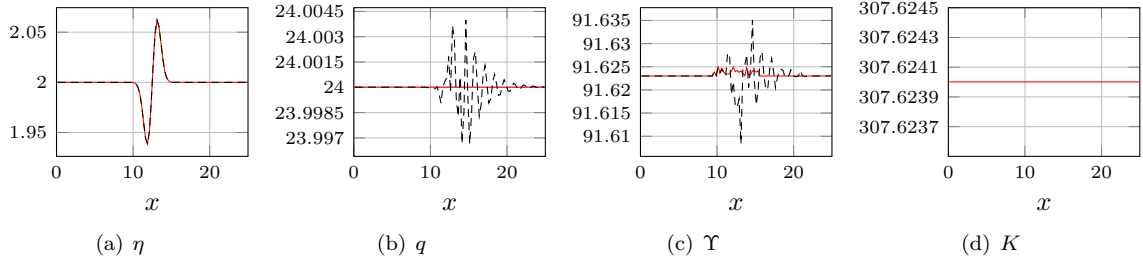
\begin{figure}
	\centering
	\subfigure[$\eta$]{
		\begin{tikzpicture}
			\begin{axis}[
				xmin=0,xmax=25,
				grid=major,
				xlabel={$x$},
				ylabel={},
				xlabel shift = 1 pt,
				ylabel shift = 1 pt,
				legend pos= south east,
				legend style={nodes={scale=0.6, transform shape}},
				tick label style={font=\scriptsize},
				width=.25\textwidth
				]
				\addplot[red]   table [y=eta   , x=x]{SUPweno5GF.dat};
				\addplot[black,densely dashed] table [y expr=\thisrow{h}+\thisrow{b}   , x=x]{SUPweno5noGF.dat};
			\end{axis}
	\end{tikzpicture}}
	\subfigure[$q$]{
		\begin{tikzpicture}
			\begin{axis}[
				xmin=0,xmax=25,
				grid=major,
				xlabel={$x$},
				ylabel={},
				xlabel shift = 1 pt,
				ylabel shift = 1 pt,
				tick label style={font=\scriptsize},
				ytick distance=0.0015,
				y tick label style={/pgf/number format/.cd,fixed relative,precision=6},
				legend pos= south east,
				legend style={nodes={scale=0.6, transform shape}},
				width=.25\textwidth
				]
				\addplot[red]   table [y=q   , x=x]{SUPweno5GF.dat};
				\addplot[black,densely dashed] table [y expr=\thisrow{h}*\thisrow{u}   , x=x]{SUPweno5noGF.dat};
			\end{axis}
	\end{tikzpicture}}
	\subfigure[$\Upsilon$]{
		\begin{tikzpicture}
			\begin{axis}[
				xmin=0,xmax=25,
				grid=major,
				xlabel={$x$},
				ylabel={},
				xlabel shift = 1 pt,
				ylabel shift = 1 pt,
				tick label style={font=\scriptsize},
				ytick distance=0.005,
				y tick label style={/pgf/number format/.cd,fixed relative,precision=6},
				legend pos= south east,
				legend style={nodes={scale=0.6, transform shape}},
				width=.25\textwidth
				]
				\addplot[red]   table [y expr=0.5*\thisrow{q}^2/\thisrow{h}^2+9.812*(\thisrow{h}+\thisrow{b})   , x=x]{SUPweno5GF.dat};
				\addplot[black,densely dashed] table [y expr=\thisrow{u}^2/2+9.812*(\thisrow{h}+\thisrow{b})   , x=x]{SUPweno5noGF.dat};
			\end{axis}
	\end{tikzpicture}}
	\subfigure[$K$]{
		\begin{tikzpicture}
			\begin{axis}[
				xmin=0,xmax=25,
				ymin=307.6235,ymax=307.6245,
				grid=major,
				xlabel={$x$},
				ylabel={},
				xlabel shift = 1 pt,
				ylabel shift = 1 pt,
				tick label style={font=\scriptsize},
				ytick distance=0.0002,
				y tick label style={/pgf/number format/.cd,fixed relative,precision=7},
				legend pos= south east,
				legend style={nodes={scale=0.6, transform shape}},
				width=.25\textwidth
				]
				\addplot[red]   table [y=k, x=x]{SUPweno5GF.dat};
			\end{axis}
	\end{tikzpicture}}
	\caption{Supercritical flow: relevant variables computed with GF-WENO5 (red continuous line) and WENO5 (black dashed line) schemes with $N_e=100$.}\label{SUP: solution}
\end{figure}
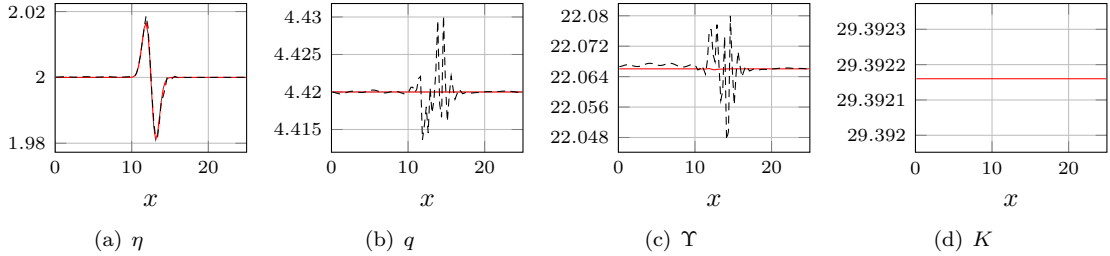
\begin{figure}
\centering
\subfigure[$\eta$]{
	\begin{tikzpicture}
		\begin{axis}[
			xmin=0,xmax=25,
			grid=major,
			xlabel={$x$},
			ylabel={},
			xlabel shift = 1 pt,
			ylabel shift = 1 pt,
			legend pos= south east,
			legend style={nodes={scale=0.6, transform shape}},
			tick label style={font=\scriptsize},
			width=.25\textwidth
			]
			\addplot[red]   table [y=eta   , x=x]{SUBweno5GF.dat};
			\addplot[black,densely dashed] table [y expr=\thisrow{h}+\thisrow{b}   , x=x]{SUBweno5noGF.dat};
		\end{axis}
\end{tikzpicture}}
\subfigure[$q$]{
	\begin{tikzpicture}
		\begin{axis}[
			xmin=0,xmax=25,
			grid=major,
			xlabel={$x$},
			ylabel={},
			xlabel shift = 1 pt,
			ylabel shift = 1 pt,
			tick label style={font=\scriptsize},
			ytick distance=0.005,
			y tick label style={/pgf/number format/.cd,fixed relative,precision=6},
			legend pos= south east,
			legend style={nodes={scale=0.6, transform shape}},
			width=.25\textwidth
			]
			\addplot[red]   table [y=q   , x=x]{SUBweno5GF.dat};
			\addplot[black,densely dashed] table [y expr=\thisrow{h}*\thisrow{u}   , x=x]{SUBweno5noGF.dat};
		\end{axis}
\end{tikzpicture}}
	\subfigure[$\Upsilon$]{
		\begin{tikzpicture}
			\begin{axis}[
				xmin=0,xmax=25,
				grid=major,
				xlabel={$x$},
				ylabel={},
				xlabel shift = 1 pt,
				ylabel shift = 1 pt,
				tick label style={font=\scriptsize},
				ytick distance=0.008,
				y tick label style={/pgf/number format/.cd,fixed relative,precision=6},
				legend pos= south east,
				legend style={nodes={scale=0.6, transform shape}},
				width=.25\textwidth
				]
				\addplot[red]   table [y expr=0.5*\thisrow{q}^2/\thisrow{h}^2+9.812*(\thisrow{h}+\thisrow{b})   , x=x]{SUBweno5GF.dat};
				\addplot[black,densely dashed] table [y expr=\thisrow{u}^2/2+9.812*(\thisrow{h}+\thisrow{b})   , x=x]{SUBweno5noGF.dat};
			\end{axis}
	\end{tikzpicture}}
	\subfigure[$K$]{
		\begin{tikzpicture}
			\begin{axis}[
				xmin=0,xmax=25,
				ymin=29.392,ymax=29.3924,
				grid=major,
				xlabel={$x$},
				ylabel={},
				xlabel shift = 1 pt,
				ylabel shift = 1 pt,
				tick label style={font=\scriptsize},
				ytick distance=0.000095,
				y tick label style={/pgf/number format/.cd,fixed relative,precision=6},
				legend pos= south east,
				legend style={nodes={scale=0.6, transform shape}},
				width=.25\textwidth
				]
				\addplot[red]   table [y=k, x=x]{SUBweno5GF.dat};
			\end{axis}
	\end{tikzpicture}}
	\caption{Subcritical flow: relevant variables computed with GF-WENO5 (red continuous line) and WENO5 (black dashed line) schemes with $N_e=100$.}\label{SUB: solution}
\end{figure}

\subsection{Perturbation of steady states without friction ($n=0$)}\label{se:pertSteadyNoFriction}

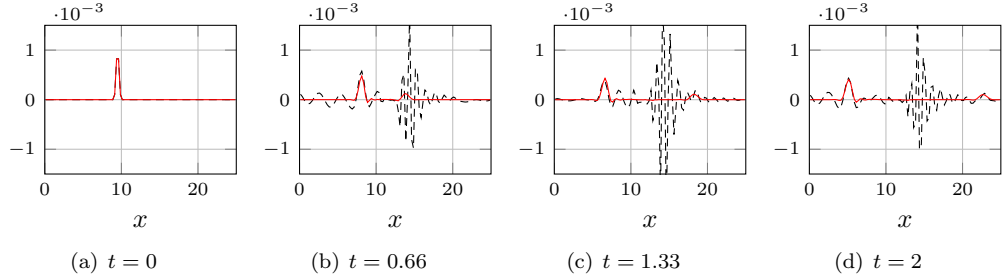
\begin{figure}
	\newcommand{\testNamenoGF}{PERTSUBweno5noGF}
	\newcommand{\testNameGF}{PERTSUBweno5GF}
	\newcommand{\Nmesh}{N100}
	\newcommand{\scale}{1.5e-3}
	\centering
	\subfigure[$t=0$]{
		\begin{tikzpicture}
			\begin{axis}[
				ymin=-\scale,ymax=\scale,
				xmin=0,xmax=25,
				grid=major,
				xlabel={$x$},
				ylabel={},
				xlabel shift = 1 pt,
				ylabel shift = 1 pt,
				legend pos= south east,
				legend style={nodes={scale=0.6, transform shape}},
				tick label style={font=\scriptsize},
				width=.25\textwidth
				]
				\addplot[black,densely dashed]   table [y=pert, x=x]{\testNamenoGF_\Nmesh_0.dat};
				\addplot[red]   table [y=pert, x=x]{\testNameGF_\Nmesh_0.dat};
			\end{axis}
	\end{tikzpicture}}
	\subfigure[$t=0.66$]{
		\begin{tikzpicture}
			\begin{axis}[
				ymin=-\scale,ymax=\scale,
				xmin=0,xmax=25,
				grid=major,
				xlabel={$x$},
				ylabel={},
				xlabel shift = 1 pt,
				ylabel shift = 1 pt,
				legend pos= south east,
				legend style={nodes={scale=0.6, transform shape}},
				tick label style={font=\scriptsize},
				width=.25\textwidth
				]
				\addplot[black,densely dashed]   table [y=pert, x=x]{\testNamenoGF_\Nmesh_1.dat};
				\addplot[red]   table [y=pert, x=x]{\testNameGF_\Nmesh_1.dat};
			\end{axis}
	\end{tikzpicture}}
	\subfigure[$t=1.33$]{
		\begin{tikzpicture}
			\begin{axis}[
				ymin=-\scale,ymax=\scale,
				xmin=0,xmax=25,
				grid=major,
				xlabel={$x$},
				ylabel={},
				xlabel shift = 1 pt,
				ylabel shift = 1 pt,
				legend pos= south east,
				legend style={nodes={scale=0.6, transform shape}},
				tick label style={font=\scriptsize},
				width=.25\textwidth
				]
				\addplot[black,densely dashed]   table [y=pert, x=x]{\testNamenoGF_\Nmesh_2.dat};
				\addplot[red]   table [y=pert, x=x]{\testNameGF_\Nmesh_2.dat};
			\end{axis}
	\end{tikzpicture}}
	\subfigure[$t=2$]{
		\begin{tikzpicture}
			\begin{axis}[
				ymin=-\scale,ymax=\scale,
				xmin=0,xmax=25,
				grid=major,
				xlabel={$x$},
				ylabel={},
				xlabel shift = 1 pt,
				ylabel shift = 1 pt,
				legend pos= south east,
				legend style={nodes={scale=0.6, transform shape}},
				tick label style={font=\scriptsize},
				width=.25\textwidth
				]
				\addplot[black,densely dashed]   table [y=pert, x=x]{\testNamenoGF_\Nmesh_3.dat};
				\addplot[red]   table [y=pert, x=x]{\testNameGF_\Nmesh_3.dat};
			\end{axis}
	\end{tikzpicture}}
	\caption{Small perturbation of the subcritical solution computed with the WENO5 scheme (black dashed) and GF WENO5 (red continuous): $h-h_{eq}$ with $N_e=100$.}\label{SUB perturbation:weno5} 
\end{figure}

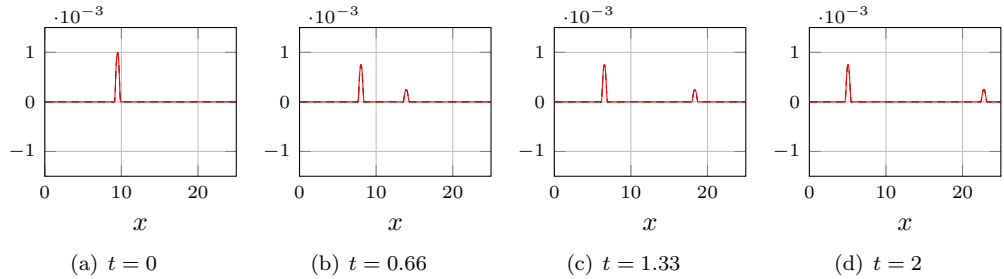
\begin{figure}
	\newcommand{\testNamenoGF}{PERTSUBweno5noGF}
	\newcommand{\testNameGF}{PERTSUBweno5GF}
	\newcommand{\Nmesh}{N800}
	\newcommand{\scale}{1.5e-3}
	\centering
	\subfigure[$t=0$]{
		\begin{tikzpicture}
			\begin{axis}[
				ymin=-\scale,ymax=\scale,
				xmin=0,xmax=25,
				grid=major,
				xlabel={$x$},
				ylabel={},
				xlabel shift = 1 pt,
				ylabel shift = 1 pt,
				legend pos= south east,
				legend style={nodes={scale=0.6, transform shape}},
				tick label style={font=\scriptsize},
				width=.25\textwidth
				]
				\addplot[black,densely dashed]   table [y=pert, x=x]{\testNamenoGF_\Nmesh_0.dat};
				\addplot[red]   table [y=pert, x=x]{\testNameGF_\Nmesh_0.dat};
			\end{axis}
	\end{tikzpicture}}
	\subfigure[$t=0.66$]{
		\begin{tikzpicture}
			\begin{axis}[
				ymin=-\scale,ymax=\scale,
				xmin=0,xmax=25,
				grid=major,
				xlabel={$x$},
				ylabel={},
				xlabel shift = 1 pt,
				ylabel shift = 1 pt,
				legend pos= south east,
				legend style={nodes={scale=0.6, transform shape}},
				tick label style={font=\scriptsize},
				width=.25\textwidth
				]
				\addplot[black,densely dashed]   table [y=pert, x=x]{\testNamenoGF_\Nmesh_1.dat};
				\addplot[red]   table [y=pert, x=x]{\testNameGF_\Nmesh_1.dat};
			\end{axis}
	\end{tikzpicture}}
	\subfigure[$t=1.33$]{
		\begin{tikzpicture}
			\begin{axis}[
				ymin=-\scale,ymax=\scale,
				xmin=0,xmax=25,
				grid=major,
				xlabel={$x$},
				ylabel={},
				xlabel shift = 1 pt,
				ylabel shift = 1 pt,
				legend pos= south east,
				legend style={nodes={scale=0.6, transform shape}},
				tick label style={font=\scriptsize},
				width=.25\textwidth
				]
				\addplot[black,densely dashed]   table [y=pert, x=x]{\testNamenoGF_\Nmesh_2.dat};
				\addplot[red]   table [y=pert, x=x]{\testNameGF_\Nmesh_2.dat};
			\end{axis}
	\end{tikzpicture}}
	\subfigure[$t=2$]{
		\begin{tikzpicture}
			\begin{axis}[
				ymin=-\scale,ymax=\scale,
				xmin=0,xmax=25,
				grid=major,
				xlabel={$x$},
				ylabel={},
				xlabel shift = 1 pt,
				ylabel shift = 1 pt,
				legend pos= south east,
				legend style={nodes={scale=0.6, transform shape}},
				tick label style={font=\scriptsize},
				width=.25\textwidth
				]
				\addplot[black,densely dashed]   table [y=pert, x=x]{\testNamenoGF_\Nmesh_3.dat};
				\addplot[red]   table [y=pert, x=x]{\testNameGF_\Nmesh_3.dat};
			\end{axis}
	\end{tikzpicture}}
	\caption{Small perturbation of the subcritical solution computed with the WENO5 scheme (black dashed) and GF WENO5 (red continuous): $h-h_{eq}$ with $N_e=800$.}\label{SUB perturbation_N800:weno5} 
\end{figure}

In this section, we add a perturbation to the tests of \cref{se:steadyNoFriction} to compare the GF-WENO5 and the classical WENO5 methods for supercritical and subcritical flows.
For the subcritical case we use the perturbation \eqref{eq:perturbation} with $\alpha=10^{-3}$. In \cref{SUB perturbation:weno5} we plot the solution for both methods at different timesteps with $N_e=100$, while in \cref{SUB perturbation_N800:weno5} we use $N_e=800$. It can be noticed that both methods converge towards the exact solution, but with the GF-WENO5 method, even with a coarse mesh we obtain a very accurate approximation of the perturbation, while the WENO5 method performs poorly. A much more refined mesh is needed for the classical methods to nicely approximate this kind of solutions.\\

\begin{figure}
	\newcommand{\testNamenoGF}{PERTSUPERweno5noGF}
	\newcommand{\testNameGF}{PERTSUPERweno5GF}
	\newcommand{\Nmesh}{N800}
	\newcommand{\scale}{2e-4}
	\centering
	\subfigure[$t=0$]{
		\begin{tikzpicture}
			\begin{axis}[
				ymin=-\scale,ymax=\scale,
				xmin=0,xmax=25,
				grid=major,
				xlabel={$x$},
				ylabel={},
				xlabel shift = 1 pt,
				ylabel shift = 1 pt,
				legend pos= south east,
				legend style={nodes={scale=0.6, transform shape}},
				tick label style={font=\scriptsize},
				width=.25\textwidth
				]
				\addplot[black,densely dashed]   table [y=pert, x=x]{\testNamenoGF_\Nmesh_0.dat};
				\addplot[red]   table [y=pert, x=x]{\testNameGF_\Nmesh_0.dat};
			\end{axis}
	\end{tikzpicture}}
	\subfigure[$t=0.33$]{
		\begin{tikzpicture}
			\begin{axis}[
				ymin=-\scale,ymax=\scale,
				xmin=0,xmax=25,
				grid=major,
				xlabel={$x$},
				ylabel={},
				xlabel shift = 1 pt,
				ylabel shift = 1 pt,
				legend pos= south east,
				legend style={nodes={scale=0.6, transform shape}},
				tick label style={font=\scriptsize},
				width=.25\textwidth
				]
				\addplot[black,densely dashed]   table [y=pert, x=x]{\testNamenoGF_\Nmesh_1.dat};
				\addplot[red]   table [y=pert, x=x]{\testNameGF_\Nmesh_1.dat};
			\end{axis}
	\end{tikzpicture}}
	\subfigure[$t=0.66$]{
		\begin{tikzpicture}
			\begin{axis}[
				ymin=-\scale,ymax=\scale,
				xmin=0,xmax=25,
				grid=major,
				xlabel={$x$},
				ylabel={},
				xlabel shift = 1 pt,
				ylabel shift = 1 pt,
				legend pos= south east,
				legend style={nodes={scale=0.6, transform shape}},
				tick label style={font=\scriptsize},
				width=.25\textwidth
				]
				\addplot[black,densely dashed]   table [y=pert, x=x]{\testNamenoGF_\Nmesh_2.dat};
				\addplot[red]   table [y=pert, x=x]{\testNameGF_\Nmesh_2.dat};
			\end{axis}
	\end{tikzpicture}}
	\subfigure[$t=1$]{
		\begin{tikzpicture}
			\begin{axis}[
				ymin=-\scale,ymax=\scale,
				xmin=0,xmax=25,
				grid=major,
				xlabel={$x$},
				ylabel={},
				xlabel shift = 1 pt,
				ylabel shift = 1 pt,
				legend pos= south east,
				legend style={nodes={scale=0.6, transform shape}},
				tick label style={font=\scriptsize},
				width=.25\textwidth
				]
				\addplot[black,densely dashed]   table [y=pert, x=x]{\testNamenoGF_\Nmesh_3.dat};
				\addplot[red]   table [y=pert, x=x]{\testNameGF_\Nmesh_3.dat};
			\end{axis}
	\end{tikzpicture}}
	\caption{Small perturbation of the supercritical solution computed with the WENO5 scheme (black dashed) and GF WENO5 (red continuous): $h-h_{eq}$ with $N_e=800$.}\label{SUPER perturbation_N800:weno5} 
\end{figure}
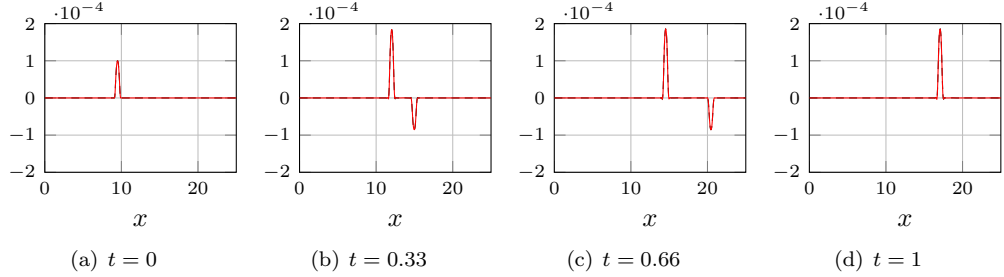

For the supercritical case we use the perturbation \eqref{eq:perturbation} with $\alpha=10^{-4}$. 
In \cref{SUPER perturbation:weno5} we plot the solution for both methods at different timesteps with $N_e=100$, while in \cref{SUPER perturbation_N800:weno5} we use $N_e=800$. 
Similar conclusions can be drawn from these simulations: the GF-WENO5 is accurate in representing the perturbation even for very coarse meshes, while the WENO5 needs more discretization cells.

\begin{figure}
	\newcommand{\testNamenoGF}{PERTSUPERweno5noGF}
	\newcommand{\testNameGF}{PERTSUPERweno5GF}
	\newcommand{\Nmesh}{N100}
	\newcommand{\scale}{2e-4}
	\centering
	\subfigure[$t=0$]{
		\begin{tikzpicture}
			\begin{axis}[
				ymin=-\scale,ymax=\scale,
				xmin=0,xmax=25,
				grid=major,
				xlabel={$x$},
				ylabel={},
				xlabel shift = 1 pt,
				ylabel shift = 1 pt,
				legend pos= south east,
				legend style={nodes={scale=0.6, transform shape}},
				tick label style={font=\scriptsize},
				width=.25\textwidth
				]
				\addplot[black,densely dashed]   table [y=pert, x=x]{\testNamenoGF_\Nmesh_0.dat};
				\addplot[red]   table [y=pert, x=x]{\testNameGF_\Nmesh_0.dat};
			\end{axis}
	\end{tikzpicture}}
	\subfigure[$t=0.33$]{
		\begin{tikzpicture}
			\begin{axis}[
				ymin=-\scale,ymax=\scale,
				xmin=0,xmax=25,
				grid=major,
				xlabel={$x$},
				ylabel={},
				xlabel shift = 1 pt,
				ylabel shift = 1 pt,
				legend pos= south east,
				legend style={nodes={scale=0.6, transform shape}},
				tick label style={font=\scriptsize},
				width=.25\textwidth
				]
				\addplot[black,densely dashed]   table [y=pert, x=x]{\testNamenoGF_\Nmesh_1.dat};
				\addplot[red]   table [y=pert, x=x]{\testNameGF_\Nmesh_1.dat};
			\end{axis}
	\end{tikzpicture}}
	\subfigure[$t=0.66$]{
		\begin{tikzpicture}
			\begin{axis}[
				ymin=-\scale,ymax=\scale,
				xmin=0,xmax=25,
				grid=major,
				xlabel={$x$},
				ylabel={},
				xlabel shift = 1 pt,
				ylabel shift = 1 pt,
				legend pos= south east,
				legend style={nodes={scale=0.6, transform shape}},
				tick label style={font=\scriptsize},
				width=.25\textwidth
				]
				\addplot[black,densely dashed]   table [y=pert, x=x]{\testNamenoGF_\Nmesh_2.dat};
				\addplot[red]   table [y=pert, x=x]{\testNameGF_\Nmesh_2.dat};
			\end{axis}
	\end{tikzpicture}}
	\subfigure[$t=1$]{
		\begin{tikzpicture}
			\begin{axis}[
				ymin=-\scale,ymax=\scale,
				xmin=0,xmax=25,
				grid=major,
				xlabel={$x$},
				ylabel={},
				xlabel shift = 1 pt,
				ylabel shift = 1 pt,
				legend pos= south east,
				legend style={nodes={scale=0.6, transform shape}},
				tick label style={font=\scriptsize},
				width=.25\textwidth
				]
				\addplot[black,densely dashed]   table [y=pert, x=x]{\testNamenoGF_\Nmesh_3.dat};
				\addplot[red]   table [y=pert, x=x]{\testNameGF_\Nmesh_3.dat};
			\end{axis}
	\end{tikzpicture}}
	\caption{Small perturbation of the supercritical solution computed with the WENO5 scheme (black dashed) and GF WENO5 (red continuous): $h-h_{eq}$ with $N_e=100$.}\label{SUPER perturbation:weno5} 
\end{figure}
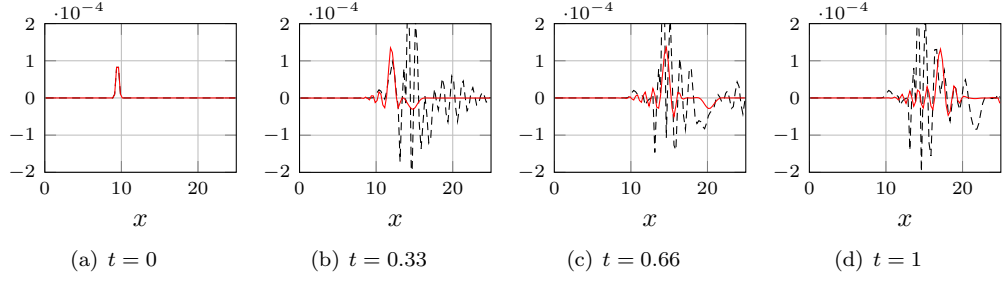

\subsection{Steady states with discontinuous bathymetry without friction ($n=0$)}\label{sec:steady_disc}
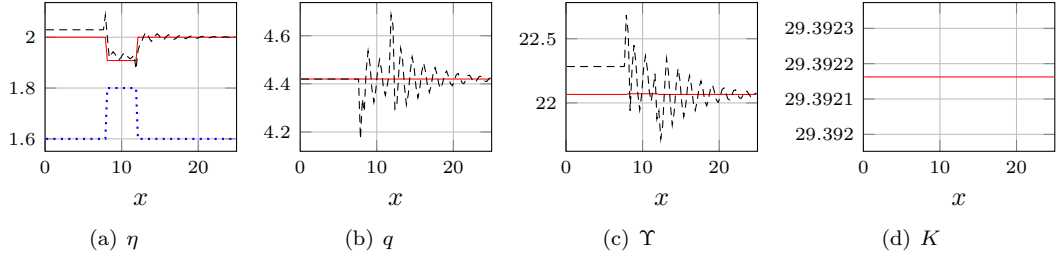
\begin{figure}
	\centering
	\newcommand{\testName}{SUB_disc_weno5}
	\newcommand{\Nmesh}{N100}
	\subfigure[$\eta$]{
		\begin{tikzpicture}
			\begin{axis}[
				xmin=0,xmax=25,
				grid=major,
				xlabel={$x$},
				ylabel={},
				xlabel shift = 1 pt,
				ylabel shift = 1 pt,
				legend pos= south east,
				legend style={nodes={scale=0.6, transform shape}},
				tick label style={font=\scriptsize},
				width=.25\textwidth
				]
				\addplot[red]   table [y=eta   , x=x]{\testName_GF_\Nmesh.dat};
				\addplot[thick,blue, dotted]   table [y expr=\thisrow{b}+1.6   , x=x]{\testName_GF_\Nmesh.dat};
				\addplot[black,densely dashed] table [y expr=\thisrow{h}+\thisrow{b}   , x=x]{\testName_noGF_\Nmesh.dat};
			\end{axis}
	\end{tikzpicture}}
\subfigure[$q$]{
	\begin{tikzpicture}
		\begin{axis}[
			xmin=0,xmax=25,
			grid=major,
			xlabel={$x$},
			ylabel={},
			xlabel shift = 1 pt,
			ylabel shift = 1 pt,
			tick label style={font=\scriptsize},
			legend pos= south east,
			legend style={nodes={scale=0.6, transform shape}},
			width=.25\textwidth
			]
			\addplot[red]   table [y=q   , x=x]{\testName_GF_\Nmesh.dat};
			\addplot[black,densely dashed] table [y expr=\thisrow{h}*\thisrow{u}   , x=x]{\testName_noGF_\Nmesh.dat};
		\end{axis}
\end{tikzpicture}}
	\subfigure[$\Upsilon$]{
		\begin{tikzpicture}
			\begin{axis}[
				xmin=0,xmax=25,
				grid=major,
				xlabel={$x$},
				ylabel={},
				xlabel shift = 1 pt,
				ylabel shift = 1 pt,
				tick label style={font=\scriptsize},
				legend pos= south east,
				legend style={nodes={scale=0.6, transform shape}},
				width=.25\textwidth
				]
				\addplot[red]   table [y expr=0.5*\thisrow{q}^2/\thisrow{h}^2+9.812*(\thisrow{h}+\thisrow{b})   , x=x]{\testName_GF_\Nmesh.dat};
				\addplot[black,densely dashed] table [y expr=\thisrow{u}^2/2+9.812*(\thisrow{h}+\thisrow{b})   , x=x]{\testName_noGF_\Nmesh.dat};
			\end{axis}
	\end{tikzpicture}}
\subfigure[$K$]{
	\begin{tikzpicture}
		\begin{axis}[
			xmin=0,xmax=25,
			ymin=29.392,ymax=29.3924,
			grid=major,
			xlabel={$x$},
			ylabel={},
			xlabel shift = 1 pt,
			ylabel shift = 1 pt,
			tick label style={font=\scriptsize},
			ytick distance=0.000095,
			y tick label style={/pgf/number format/.cd,fixed relative,precision=6},
			legend pos= south east,
			legend style={nodes={scale=0.6, transform shape}},
			width=.25\textwidth
			]
			\addplot[red]   table [y=k, x=x]{\testName_GF_\Nmesh.dat};
		\end{axis}
\end{tikzpicture}}
\caption{Subcritical flow: relevant variables computed with GF-WENO5 (red continuous line), WENO5 (black dashed line) schemes and rescaled $b$ (blue dotted line) with $N_e=100$.}\label{SUB disc: solution}
\end{figure}

Here, we change the bathymetry and test the subcritical and supercritical flows with a discontinuous step. The bathymetry is defined by
\begin{equation}
	b=\begin{cases}
		0.2 & \text{if } 8<x<12,\\
		0. &\text{else.}
	\end{cases}
\end{equation}
In \cref{SUB disc: solution} we test the subcritical flow with the same initial conditions of \cref{se:steadyNoFriction} with $T=500$ and $N_e=100$. We can observe that the GF-WENO5 performs amazingly without producing any oscillations and obtaining errors of the order of $10^{-10}$, while the WENO5 scheme wildly oscillates after the discontinuities of the bathymetry.\\

\begin{figure}
\centering
\newcommand{\testName}{SUPER_disc_weno5}
\newcommand{\Nmesh}{N100}
\subfigure[$\eta$]{
	\begin{tikzpicture}
		\begin{axis}[
			xmin=0,xmax=25,
			grid=major,
			xlabel={$x$},
			ylabel={},
			xlabel shift = 1 pt,
			ylabel shift = 1 pt,
			legend pos= south east,
			legend style={nodes={scale=0.6, transform shape}},
			tick label style={font=\scriptsize},
			width=.25\textwidth
			]
			\addplot[red]   table [y=eta   , x=x]{\testName_GF_\Nmesh.dat};
			\addplot[thick,blue, dotted]   table [y expr=\thisrow{b}+1.75   , x=x]{\testName_GF_\Nmesh.dat};
			\addplot[black,densely dashed] table [y expr=\thisrow{h}+\thisrow{b}   , x=x]{\testName_noGF_\Nmesh.dat};
		\end{axis}
\end{tikzpicture}}
\subfigure[$q$]{
	\begin{tikzpicture}
		\begin{axis}[
			xmin=0,xmax=25,
			grid=major,
			xlabel={$x$},
			ylabel={},
			xlabel shift = 1 pt,
			ylabel shift = 1 pt,
			tick label style={font=\scriptsize},
			legend pos= south east,
			legend style={nodes={scale=0.6, transform shape}},
			width=.25\textwidth
			]
			\addplot[red]   table [y=q   , x=x]{\testName_GF_\Nmesh.dat};
			\addplot[black,densely dashed] table [y expr=\thisrow{h}*\thisrow{u}   , x=x]{\testName_noGF_\Nmesh.dat};
		\end{axis}
\end{tikzpicture}}
\subfigure[$\Upsilon$]{
\begin{tikzpicture}
	\begin{axis}[
		xmin=0,xmax=25,
		grid=major,
		xlabel={$x$},
		ylabel={},
		xlabel shift = 1 pt,
		ylabel shift = 1 pt,
		tick label style={font=\scriptsize},
		legend pos= south east,
		legend style={nodes={scale=0.6, transform shape}},
		width=.25\textwidth
		]
		\addplot[red]   table [y expr=0.5*\thisrow{q}^2/\thisrow{h}^2+9.812*(\thisrow{h}+\thisrow{b})   , x=x]{\testName_GF_\Nmesh.dat};
		\addplot[black,densely dashed] table [y expr=\thisrow{u}^2/2+9.812*(\thisrow{h}+\thisrow{b})   , x=x]{\testName_noGF_\Nmesh.dat};
	\end{axis}
\end{tikzpicture}}
\subfigure[$K$]{
	\begin{tikzpicture}
		\begin{axis}[
			xmin=0,xmax=25,
			ymin=307.6235,ymax=307.6245,
			grid=major,
			xlabel={$x$},
			ylabel={},
			xlabel shift = 1 pt,
			ylabel shift = 1 pt,
			tick label style={font=\scriptsize},
			ytick distance=0.0003,
			y tick label style={/pgf/number format/.cd,fixed relative,precision=7},
			legend pos= south east,
			legend style={nodes={scale=0.6, transform shape}},
			width=.25\textwidth
			]
			\addplot[red]   table [y=k, x=x]{\testName_GF_\Nmesh.dat};
		\end{axis}
\end{tikzpicture}}
\caption{Supercritical flow: relevant variables computed with GF-WENO5 (red continuous line) and WENO5 (black dashed line) schemes and rescaled $b$ (blue dotted line)  with $N_e=100$.}\label{SUPER disc: solution}
\end{figure}
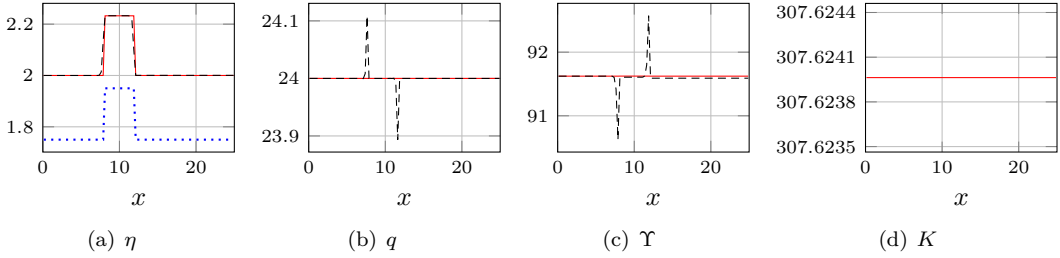

For the supercritical case we have used the initial conditions of \cref{se:steadyNoFriction} till final time $T=50$ with $N_e=100$. In \cref{SUPER disc: solution} we observe that WENO5 performs better than before, still producing spurious oscillations and not exactly catching the outflow solution. On the other side, GF-WENO5 obtains constant global fluxes with a machine precision accuracy.

\subsection{Steady states with friction ($n=0.05$)}\label{se:steadyFriction}

In this section we focus on the supercritical and subcritical flows studied in \cref{se:steadyNoFriction} including the friction term in the source term.
As for the previous cases, also when friction with constant Manning coefficient $n$ is present, we can obtain moving equilibria. 
Again the quantity that are preserved at equilibrium are $q$ and $K$.  
That is why it is interesting to perform simulations similar to the previous ones comparing standard methods with GF ones.
We consider the subcritical case defined in \eqref{eq:sub_ic} and the supercritical case defined in \eqref{eq:super_ic} with the same bathymetry, on which we add the friction term with Manning coefficient $n=0.05$. In the supercritical case, the friction term implies a slow down of the physical speed from left to right and a consequent increasing of $\eta$ from left to right. In the subcritical case, conversely, we expect $h$ to decrease from left to right and the speed to increase. The variable $\Upsilon$ is not conserved  and, for these tests, there is not another constant variable that can be easily computed analytically, hence, we do not plot it. \\

\begin{figure}
	\centering
	\subfigure[$\eta$]{
		\begin{tikzpicture}
			\begin{axis}[
				xmin=0,xmax=25,
				grid=major,
				xlabel={$x$},
				ylabel={},
				xlabel shift = 1 pt,
				ylabel shift = 1 pt,
				legend pos= south east,
				legend style={nodes={scale=0.6, transform shape}},
				tick label style={font=\scriptsize},
				width=.32\textwidth
				]
				\addplot[red]   table [y=eta   , x=x]{SUPFRICweno5GF.dat};
				\addplot[black,densely dashed] table [y=eta   , x=x]{SUPFRICweno5noGF.dat};
			\end{axis}
	\end{tikzpicture}}
	\subfigure[$q$]{
		\begin{tikzpicture}
			\begin{axis}[
				xmin=0,xmax=25,
				grid=major,
				xlabel={$x$},
				ylabel={},
				xlabel shift = 1 pt,
				ylabel shift = 1 pt,
				tick label style={font=\scriptsize},
				legend pos= south east,
				legend style={nodes={scale=0.6, transform shape}},
				width=.32\textwidth
				]
				\addplot[red]   table [y=q   , x=x]{SUPFRICweno5GF.dat};
				\addplot[black,densely dashed] table [y expr=\thisrow{h}*\thisrow{u}   , x=x]{SUPFRICweno5noGF.dat};
			\end{axis}
	\end{tikzpicture}}
		\subfigure[$K$]{
			\begin{tikzpicture}
				\begin{axis}[
					xmin=0,xmax=25,
					ymin=307.4,ymax=307.8,
					grid=major,
					xlabel={$x$},
					ylabel={},
					xlabel shift = 1 pt,
					ylabel shift = 1 pt,
					tick label style={font=\scriptsize},
					legend pos= south east,
					legend style={nodes={scale=0.6, transform shape}},
					width=.32\textwidth
					]
					\addplot[red]   table [y=k, x=x]{SUPFRICweno5GF.dat};
				\end{axis}
		\end{tikzpicture}}
		\caption{Supercritical flow with friction: relevant variables computed with GF-WENO5 (red continuous line) and WENO5 (black dashed line) schemes with $N_e=100$.}\label{SUPFRIC: solution}
	\end{figure}
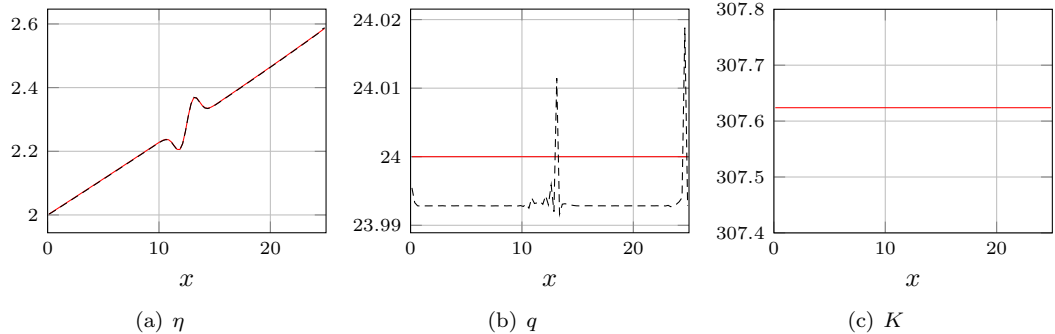

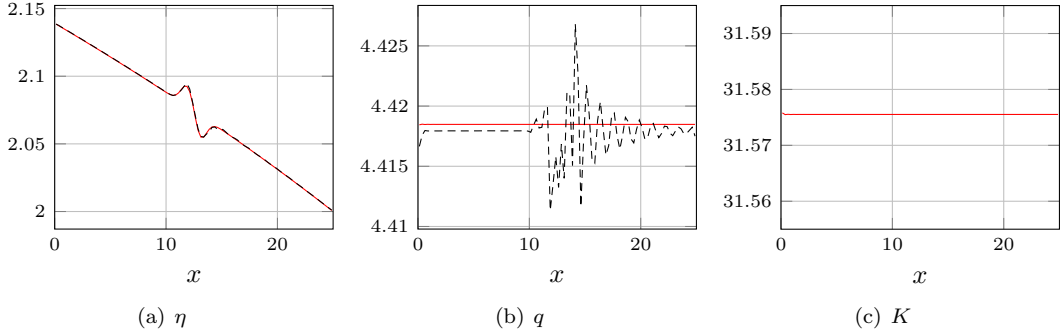
\begin{figure}
	\centering
	\subfigure[$\eta$]{
		\begin{tikzpicture}
			\begin{axis}[
				xmin=0,xmax=25,
				grid=major,
				xlabel={$x$},
				ylabel={},
				xlabel shift = 1 pt,
				ylabel shift = 1 pt,
				legend pos= south east,
				legend style={nodes={scale=0.6, transform shape}},
				tick label style={font=\scriptsize},
				width=.32\textwidth
				]
				\addplot[red]   table [y=eta   , x=x]{SUBFRICweno5GF.dat};
				\addplot[black,densely dashed] table [y expr=\thisrow{h}+\thisrow{b}   , x=x]{SUBFRICweno5noGF.dat};
			\end{axis}
	\end{tikzpicture}}
	\subfigure[$q$]{
		\begin{tikzpicture}
			\begin{axis}[
				xmin=0,xmax=25,
				grid=major,
				xlabel={$x$},
				ylabel={},
				xlabel shift = 1 pt,
				ylabel shift = 1 pt,
				tick label style={font=\scriptsize},
				y tick label style={/pgf/number format/.cd,fixed relative,precision=4},
				legend pos= south east,
				legend style={nodes={scale=0.6, transform shape}},
				width=.32\textwidth
				]
				\addplot[red]   table [y=q   , x=x]{SUBFRICweno5GF.dat};
				\addplot[black,densely dashed] table [y expr=\thisrow{h}*\thisrow{u}   , x=x]{SUBFRICweno5noGF.dat};
			\end{axis}
	\end{tikzpicture}}
\subfigure[$K$]{
	\begin{tikzpicture}
		\begin{axis}[
			xmin=0,xmax=25,
			ymin=31.555,ymax=31.595,
			grid=major,
			xlabel={$x$},
			ylabel={},
			xlabel shift = 1 pt,
			ylabel shift = 1 pt,
			tick label style={font=\scriptsize},
			legend pos= south east,
			legend style={nodes={scale=0.6, transform shape}},
			width=.32\textwidth
			]
			\addplot[red]   table [y=k, x=x]{SUBFRICweno5GF.dat};
		\end{axis}
\end{tikzpicture}}
\caption{Subcritical flow with friction: relevant variables computed with GF-WENO5 (red continuous line) and WENO5 (black dashed line) schemes with $N_e=100$.}\label{SUBFRIC: solution}
\end{figure}

We display in Figures~\ref{SUPFRIC: solution} and~\ref{SUBFRIC: solution} both the solutions computed with the new GF-WENO5 WB scheme and the
classical WENO5. It should be noticed that both schemes obtain valid and consistent results with the expected solution.
The difference between the schemes is remarkable and the global flux variables clearly highlights it. 
The WENO5 case, without the global flux, is characterized by strong spurious oscillations around the area where the effect of the bathymetry is stronger.
On the other side, the GF-WENO5 results are very precise and are able to preserve the global flux variables up to $\sim 10^{-9}$.

\subsection{Perturbation of steady states with friction ($n = 0.05$)}

Also for these friction tests, we add the perturbation \eqref{eq:perturbation} to the final solution obtained at $T=50$ for the supercritical case with $\alpha=10^{-4}$ and at $T=200$ with $\alpha=10^{-3}$ for the subcritical case. Again, we use $N_e=100$ cells of the domain. 
In Figure~\ref{SUPER perturbation_friction:weno5}, we notice that the supercritical case with friction is the only case in which the 
classical WENO scheme is able to reach the steady state (residuals approach machine precision), hence the two solutions appear
closer with respect to other cases. Probably, the friction plays a major role in this simulation by introducing some stabilization effects.
Nevertheless, the classical WENO scheme exhibits a more oscillatory and unphysical behaviour in between the perturbation waves. 
In Figure~\ref{SUB perturbation_friction:weno5}, a much better prediction of the perturbation evolution and overall solution is obtained by
the GF-WENO5 WB scheme. The classical WENO5 scheme, in this case, is not able to properly achieve the steady state, therefore producing 
a lot of numerical errors within the computational domain. 
This spoils the overall accuracy of the solution, that can be reobtained by
refining the mesh. 
\begin{figure}
	\newcommand{\testNamenoGF}{PERTSUPERFRICTIONweno5noGF}
	\newcommand{\testNameGF}{PERTSUPERFRICTIONweno5GF}
	\newcommand{\Nmesh}{N100}
	\newcommand{\scale}{1.5e-4}
	\centering
	\subfigure[$t=0$]{
		\begin{tikzpicture}
			\begin{axis}[
				ymin=-\scale,ymax=\scale,
				xmin=0,xmax=25,
				grid=major,
				xlabel={$x$},
				ylabel={},
				xlabel shift = 1 pt,
				ylabel shift = 1 pt,
				legend pos= south east,
				legend style={nodes={scale=0.6, transform shape}},
				tick label style={font=\scriptsize},
				width=.25\textwidth
				]
				\addplot[black,densely dashed]   table [y=pert, x=x]{\testNamenoGF_\Nmesh_0.dat};
				\addplot[red]   table [y=pert, x=x]{\testNameGF_\Nmesh_0.dat};
			\end{axis}
	\end{tikzpicture}}
	\subfigure[$t=0.33$]{
		\begin{tikzpicture}
			\begin{axis}[
				ymin=-\scale,ymax=\scale,
				xmin=0,xmax=25,
				grid=major,
				xlabel={$x$},
				ylabel={},
				xlabel shift = 1 pt,
				ylabel shift = 1 pt,
				legend pos= south east,
				legend style={nodes={scale=0.6, transform shape}},
				tick label style={font=\scriptsize},
				width=.25\textwidth
				]
				\addplot[black,densely dashed]   table [y=pert, x=x]{\testNamenoGF_\Nmesh_1.dat};
				\addplot[red]   table [y=pert, x=x]{\testNameGF_\Nmesh_1.dat};
			\end{axis}
	\end{tikzpicture}}
	\subfigure[$t=0.66$]{
		\begin{tikzpicture}
			\begin{axis}[
				ymin=-\scale,ymax=\scale,
				xmin=0,xmax=25,
				grid=major,
				xlabel={$x$},
				ylabel={},
				xlabel shift = 1 pt,
				ylabel shift = 1 pt,
				legend pos= south east,
				legend style={nodes={scale=0.6, transform shape}},
				tick label style={font=\scriptsize},
				width=.25\textwidth
				]
				\addplot[black,densely dashed]   table [y=pert, x=x]{\testNamenoGF_\Nmesh_2.dat};
				\addplot[red]   table [y=pert, x=x]{\testNameGF_\Nmesh_2.dat};
			\end{axis}
	\end{tikzpicture}}
	\subfigure[$t=1$]{
		\begin{tikzpicture}
			\begin{axis}[
				ymin=-\scale,ymax=\scale,
				xmin=0,xmax=25,
				grid=major,
				xlabel={$x$},
				ylabel={},
				xlabel shift = 1 pt,
				ylabel shift = 1 pt,
				legend pos= south east,
				legend style={nodes={scale=0.6, transform shape}},
				tick label style={font=\scriptsize},
				width=.25\textwidth
				]
				\addplot[black,densely dashed]   table [y=pert, x=x]{\testNamenoGF_\Nmesh_3.dat};
				\addplot[red]   table [y=pert, x=x]{\testNameGF_\Nmesh_3.dat};
			\end{axis}
	\end{tikzpicture}}
	\caption{Small perturbation of the friction supercritical solution computed with the WENO5 scheme (black dashed) and GF WENO5 (red continuous): $h-h_{eq}$ with $N_e=100$.}\label{SUPER perturbation_friction:weno5} 
\end{figure}
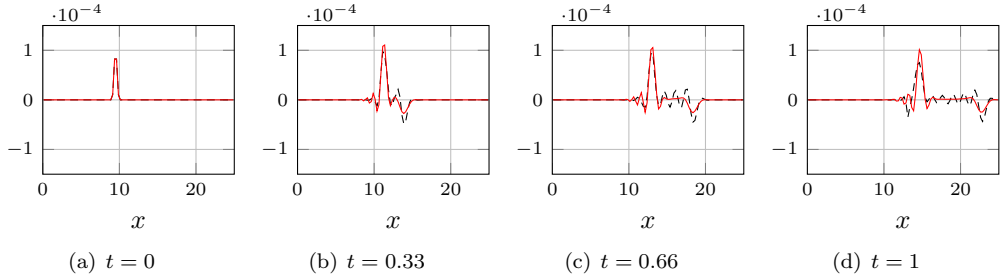
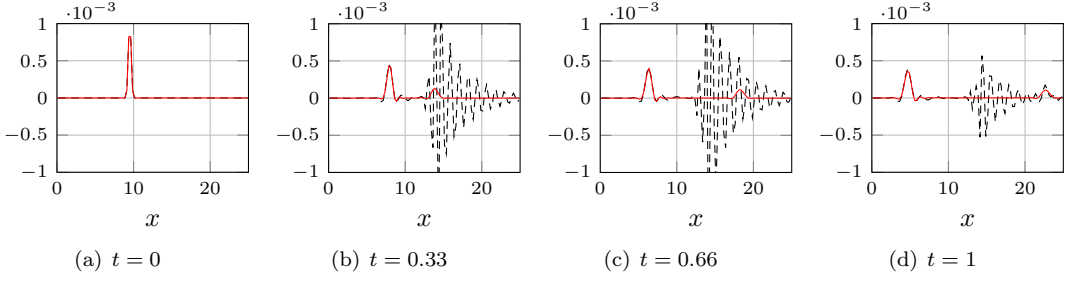
\begin{figure}
	\newcommand{\testNamenoGF}{PERTSUBFRICTIONweno5noGF}
	\newcommand{\testNameGF}{PERTSUBFRICTIONweno5GF}
	\newcommand{\Nmesh}{N100}
	\newcommand{\scale}{1e-3}
	\centering
	\subfigure[$t=0$]{
		\begin{tikzpicture}
			\begin{axis}[
				ymin=-\scale,ymax=\scale,
				xmin=0,xmax=25,
				grid=major,
				xlabel={$x$},
				ylabel={},
				xlabel shift = 1 pt,
				ylabel shift = 1 pt,
				legend pos= south east,
				legend style={nodes={scale=0.6, transform shape}},
				tick label style={font=\scriptsize},
				width=.25\textwidth
				]
				\addplot[black,densely dashed]   table [y=pert, x=x]{\testNamenoGF_\Nmesh_0.dat};
				\addplot[red]   table [y=pert, x=x]{\testNameGF_\Nmesh_0.dat};
			\end{axis}
	\end{tikzpicture}}
	\subfigure[$t=0.33$]{
		\begin{tikzpicture}
			\begin{axis}[
				ymin=-\scale,ymax=\scale,
				xmin=0,xmax=25,
				grid=major,
				xlabel={$x$},
				ylabel={},
				xlabel shift = 1 pt,
				ylabel shift = 1 pt,
				legend pos= south east,
				legend style={nodes={scale=0.6, transform shape}},
				tick label style={font=\scriptsize},
				width=.25\textwidth
				]
				\addplot[black,densely dashed]   table [y=pert, x=x]{\testNamenoGF_\Nmesh_1.dat};
				\addplot[red]   table [y=pert, x=x]{\testNameGF_\Nmesh_1.dat};
			\end{axis}
	\end{tikzpicture}}
	\subfigure[$t=0.66$]{
		\begin{tikzpicture}
			\begin{axis}[
				ymin=-\scale,ymax=\scale,
				xmin=0,xmax=25,
				grid=major,
				xlabel={$x$},
				ylabel={},
				xlabel shift = 1 pt,
				ylabel shift = 1 pt,
				legend pos= south east,
				legend style={nodes={scale=0.6, transform shape}},
				tick label style={font=\scriptsize},
				width=.25\textwidth
				]
				\addplot[black,densely dashed]   table [y=pert, x=x]{\testNamenoGF_\Nmesh_2.dat};
				\addplot[red]   table [y=pert, x=x]{\testNameGF_\Nmesh_2.dat};
			\end{axis}
	\end{tikzpicture}}
	\subfigure[$t=1$]{
		\begin{tikzpicture}
			\begin{axis}[
				ymin=-\scale,ymax=\scale,
				xmin=0,xmax=25,
				grid=major,
				xlabel={$x$},
				ylabel={},
				xlabel shift = 1 pt,
				ylabel shift = 1 pt,
				legend pos= south east,
				legend style={nodes={scale=0.6, transform shape}},
				tick label style={font=\scriptsize},
				width=.25\textwidth
				]
				\addplot[black,densely dashed]   table [y=pert, x=x]{\testNamenoGF_\Nmesh_3.dat};
				\addplot[red]   table [y=pert, x=x]{\testNameGF_\Nmesh_3.dat};
			\end{axis}
	\end{tikzpicture}}
	\caption{Small perturbation of the friction subcritical solution computed with the WENO5 scheme (black dashed) and GF WENO5 (red continuous): $h-h_{eq}$ with $N_e=100$.}\label{SUB perturbation_friction:weno5} 
\end{figure}

\section{Conclusion} \label{se:summary}

In this paper, we presented a novel arbitrary high-order well-balanced finite-volume method based on 
flux globalization for the shallow water equations with source terms.
The high order accuracy is obtained through a WENO reconstruction performed on the variables of interest, i.e.\
free surface level and global fluxes. The flux globalization allows to exactly preserve the constant fluxes in moving water equilibria.
The new scheme has been also designed to preserve the ``lake-at-rest'' solution. 
This was possible by introducing a particular quadrature procedure for the source flux and 
considering a jump of the global fluxes at each interface. 
Several tests have been performed to assess the properties of the scheme. The preservation of the constant flux has been verified on academic moving water test cases (subcritical, transcritical and supercritical), while the well-balancedness with respect to the lake-at-rest equilibrium has been tested on the classical ``lake-at-rest'' and ``perturbed lake-at-rest'' cases.
Once the space discretization is well implemented, the introduction of additional source terms is straightforward, e.g.\ Manning friction, as demonstrated in the simulation section.
With the presented scheme, we were able to outperform classical methods in many situations, obtaining much more accurate solutions and useful properties also at the discrete level. Moreover, the high order accuracy of scheme allows to obtain precise solutions also when the equilibria are not reached. Several other extensions  are being investigated. These include the DG-SEM approach of 
\cite{mantri2022} as well as a  continuous finite element formulation. The WENO work presented here could provide a setting for handling discontinuities in these works. 
Further improvements of the method might involve \textit{a posteriori} limiters to deal with wet and dry areas.
The benefits of the global flux idea in other settings (more complex equations, multi dimensional problems)  is also under investigation. 

\subsection*{Acknowledgements}
M. Ricchiuto is a member of the CARDAMOM team at INRIA and University of Bordeaux.
During the development of this project, M. Ciallella was funded by an INRIA PhD fellowship and by a postdoctoral fellowship at ENSAM (I2M).
D. Torlo has benefitted of INRIA postdoctoral fellowship, and by a postdoctoral fellowship in SISSA.

\appendix
\section{Positive reconstruction of water height at interfaces}\label{sec:positiveHreconstruction}
\begin{proof}
From the cell averages a high order WENO reconstruction is performed on the fluxes to have $q^{L,R}_{\iip}$ and $K^{L,R}_{\iip}$.
Equipped with $q^{L,R}_{\iip},K^{L,R}_{\iip},\mathcal R^{L,R}_{\iip}$, the point values $h^{L,R}_{\iip}$ can be obtained by solving the nonlinear equation coming from the 
definition of the global variable $K$ in~\eqref{eq:globalCL}:
\begin{equation}\label{eq:cubic for h}
	K^L_{\iip} = \frac{\left(q^L_{\iip}\right)^2}{h^L_{\iip}} + \frac{g}{2}\left(h^L_{\iip}\right)^2 + \mathcal R^L_{\iip}, \quad K^R_{\iip} = \frac{\left(q^R_{\iip}\right)^2}{h^R_{\iip}} + \frac{g}{2}\left(h^R_{\iip}\right)^2 + \mathcal R^R_{\iip}
\end{equation}
Let us solve the depressed cubic equation~\eqref{eq:cubic for h} for $h^L_{\iip}$ (the solution for $h^R_{\iip}$ is obtained with the same procedure).
First of all, it can be noticed that~\eqref{eq:cubic for h} does not have any positive solution unless the determinant is greater than zero, meaning that
\begin{equation}\label{eq:determinant cubic}
	\left(q^L_{\iip}\right)^4 < \frac{8 \left(K^L_{\iip}-\mathcal R^L_{\iip}\right)^3}{27\,g}.
\end{equation}
If~\eqref{eq:determinant cubic} is not satisfied we reconstruct $\eta^L_{\iip}$ and then compute $h^L_{\iip}$ given that 
\begin{equation}\label{eq:hroot}
	h^L_{\iip}=\eta^L_{\iip}-b^L_{\iip}.
\end{equation}
If~\eqref{eq:determinant cubic} is satisfied, then we have to deal with two possibilities. First, if $q^L_{\iip}=0$, we obtain the unique positive solution
$$ h^L_{\iip} = \sqrt{\frac{2\left(K^L_{\iip}-\mathcal R^L_{\iip}\right)}{g}}, $$
while if $q^L_{\iip}\neq0$, we solve Equation~\eqref{eq:cubic for h} for $h^L_{\iip}$ and obtain the following three solutions:
\begin{equation}
	h^L_{\iip} = 2\sqrt{P}\cos\left(\frac{1}{3}\left[\Theta+2\pi k\right]\right), \quad k=0,1,2,
\end{equation}
where
\begin{equation}
	P:=\frac{2\left(K^L_{\iip}-R^L_{\iip}\right)}{3g} \quad\text{ and }\quad \Theta:=\arccos\left(-\frac{\left(q^L_{\iip}\right)^2}{g\,P^{3/2}}\right).
\end{equation}
It can be shown that one of these roots is negative, whilst the other two roots, corresponding to the subcritical and supercritical cases, are positive.
We choose the one closer to the corresponding value of $h^L_{\iip}$ given in~\eqref{eq:hroot}.
\section{Proof of proposition \ref{prop:lar}}\label{app:lar}
To prove the result we will prove that,  given WENO polynomials of $h$, $q$, and $b$, for the lake at rest state $[ h(x)=\eta_0-b(x)\;\;q(x)=0]^t$ we have
$\bar{\mathcal{G}}_{i+1}=\bar{\mathcal{G}}_i$ $\forall\, i$. Clearly, the first component is the reconstruction of $q\equiv 0$, hence it is zero. Let us focus on $K$, the second component of the global flux.
It is then enough to show that $\bar K_i$ is constant to have that $\bar q_i=0$ $\forall\, i$ at the following timesteps. Let us first compute the values of $K$ at each
quadrature point. Simple computations show that
\begin{align}\label{eq:lakeatrest0}
	K_{i,q} =\mathcal{F}_{i,q} +\mathcal R_{i,q} \;=& \mathcal R^R_{\iin} + 
	g\frac{(\eta_0-b_{i,q})^2}{2} + g\eta_0\left(\tilde b_{i,q} - b_{\iin}^R\right)  - g\left(\frac{(\tilde b_{i,q})^2}{2} - \frac{(b_{\iin}^R)^2}{2}\right) \\
	=&\mathcal R^R_{\iin} +g \frac{\eta_0^2}{2} -g\eta_0  b_{\iin}^R +g \frac{(b_{\iin}^R)^2}{2}.
\end{align}
This shows that $K$ is constant across the quadrature points, and thus $\bar K_i = K_{i,q}$.
We now need to show that this constant is the same $\forall \,i$. This is readily shown by substitution of the relevant quantities and some simple algebra: 
\begin{subequations}
\begin{align}
		\bar K_{i+1}-\bar K_i \;=&\; \mathcal R^R_{\iip} - \mathcal R^R_{\iin} -g\eta_0  b_{\iip}^R +g \frac{(b_{\iip}^R)^2}{2} +g\eta_0  b_{\iin}^R -g \frac{(b_{\iin}^R)^2}{2}    \\
		=& \mathcal R^L_{\iip} - \mathcal R^R_{\iin} + [\![ \mathcal R_{\iip} ]\!] -g\eta_0  b_{\iip}^R +g \frac{(b_{\iip}^R)^2}{2} +g\eta_0  b_{\iin}^R -g \frac{(b_{\iin}^R)^2}{2}   \\	
		=& \underbrace{g\eta_0\left(b_\iip^L - b_\iin^R\right) - g\left(\frac{(b_\iip^L)^2}{2} - \frac{(b_\iin^R)^2}{2}\right)}_{\mathcal R^L_{\iip} - \mathcal R^R_{\iin}} +[\![ \mathcal R_{\iip} ]\!]\\
		-&g\eta_0  b_{\iip}^R +g \frac{(b_{\iip}^R)^2}{2} +g\eta_0  b_{\iin}^R -g \frac{(b_{\iin}^R)^2}{2}  =\\
		=& g\eta_0\left(b_\iip^L - b_\iin^R\right) - g\left(\frac{(b_\iip^L)^2}{2} - \frac{(b_\iin^R)^2}{2}\right) +[\![ \mathcal R_{\iip} ]\!]=0,\label{eq:last_defi_jump}
\end{align}
\end{subequations}
recalling that 
\begin{equation}
	[\![ \mathcal R_{\iip} ]\!]=g\eta_0\left(b_\iip^R - b_\iin^L\right) - g\left(\frac{(b_\iip^R)^2}{2} - \frac{(b_\iin^L)^2}{2}\right),
\end{equation}
which achieves the proof.
\end{proof}
\begin{remark}[Definition of the jump of $\mathcal{R}$]\label{rem:definition_jump_R}
	Clearly, the definition of $[\![ \mathcal R_{\iip} ]\!]$ is obtained with the goal of having $\bar K_{i+1}-\bar K_i=0$ for lake at rest equilibrium, and it can be easily derived from \eqref{eq:last_defi_jump}.
\end{remark}

 \bibliographystyle{abbrv}
\bibliography{literature}
\end{document}